\documentclass[12pt]{amsart}
\usepackage{amssymb,amscd,verbatim}
\swapnumbers

\renewcommand{\d}{\mbox{\LARGE $\cdot $}}
\newcommand{\Xs}{X_{\d}}            %For simplicial
\renewcommand{\hat}{\widehat}
\newcommand{\Gal}{{\rm Gal}\,}       % For Galois
\newcommand{\Lie}{{\rm Lie}\,}       % For Lie 
\newcommand{\Spec}{\operatorname{Spec}} %for Spec
\newcommand{\Hom}{\operatorname{Hom}}      % For Hom
\newcommand{\Ext}{\operatorname{Ext}}      % For Ext
\newcommand{\Biext}{\operatorname{Biext}}      % For Biext
\newcommand{\RHom}{\operatorname{RHom}}    % For RHom
\newcommand{\DM}{{\rm DM}}          % For Voevodsky DM
\newcommand{\M}{\mathcal{M}_1}   % For category of 1-Motives
\newcommand{\eM}{\mathcal{M}_1^{\rm eff}} % For effective 1-Motives
\newcommand{\1}{{}_{\leq 1}}
\newcommand{\car}{\operatorname{char}}
\newcommand{\red} {{\rm red}\,}
\newcommand{\Div}{\operatorname{Div}}

\newcommand{\sext}{\text{${\mathcal E}xt\,$}}  % For sheaf Ext
\newcommand{\stor}{\text{${\mathcal T}or\,$}}  % For sheaf Tor
\newcommand{\shom}{\text{${\mathcal H}om\,$}}  %For sheaf Hom
\newcommand{\ihom}{{\rm\underline{Hom}}}  % For internal Hom
\newcommand{\oo}{\operatornamewithlimits{\otimes}\limits}

\newcommand{\Aff}{\mathbb{A}}   % For Affine space
\newcommand{\Chi}{\mbox{\large $\chi$}} % For characters

\newcommand{\C}{\mathbb{C}}     % For Complex numbers
\newcommand{\Q}{\mathbb{Q}}     % For rational numbers
\newcommand{\Z}{\mathbb{Z}}     % For integers

\newcommand{\D}{\mathbb{D}}    %Dieudonne
\newcommand{\W}{\mathbb{W}}    %Witt vectors
\newcommand{\G}{\mathbb{G}}     % For group schemes
\newcommand{\HH}{\mathbb{H}}    % For Hyper
\newcommand{\EExt}{{\rm \mathbb{E}xt}} % For HyperExtension
\newcommand{\bPic}{{\rm\mathbb{P}ic}}    % For SimpPic

\newcommand{\im}{\operatorname{Im}}        % For the image of a morphism
\renewcommand{\ker}{\operatorname{Ker}}  % For the kernel of a morphism
\newcommand{\coker}{\operatorname{Coker}} % For the cokernel of a morphism
\newcommand{\gr}{\operatorname{gr}}        % For associated Graded
\newcommand{\Pic}{{\rm Pic}}     % For the Picard variety
\newcommand{\Alb}{{\rm Alb}}     % For the Albanese variety
\newcommand{\RPic}{{\rm RPic}}     % For the Picard complex
\newcommand{\LAlb}{{\rm LAlb}}     % For the Albanese complex
\newcommand{\LA}[1]{\mbox{${\rm L}_{#1}{\rm Alb}$}}
\newcommand{\RA}[1]{\mbox{${\rm R}^{#1}{\rm Pic}$}}

\newcommand{\Tot}{{\rm Tot}}     % For the Total Complex
\newcommand{\NS}  {{\rm NS}}      % For the N\'eron Severi group
\newcommand{\qi}{{\rm q.i.}\,}      % For quasi-iso
\newcommand{\ssp}[1]{\mbox{$\scriptscriptstyle {#1}$}}
\newcommand{\by}[1]{\stackrel{#1}{\rightarrow}}
\newcommand{\longby}[1]{\stackrel{#1}{\longrightarrow}}

\newcommand{\iso}{\stackrel{\sim}{\longrightarrow}}

\renewcommand{\tilde}{\widetilde}
\newcommand{\df}{\mbox{\,${:=}$}\,}
\newcommand{\ie}{{\it i.e.\/},\ }
\newcommand{\cf}{{\it cf.\/}\ }
\newcommand{\eg}{{\it e.g.\/},\ }
\newcommand{\et} {\mbox{\scriptsize{\rm {\'e}t}}}

\newcommand{\fr}{\mbox{\scriptsize{\rm fr}}}
\newcommand{\tor}{\mbox{\scriptsize{\rm tor}}}

\newcommand{\tf}{\mbox{\scriptsize{\rm tf}}}
\newcommand{\eff}{\mbox{\scriptsize{\rm eff}}}
\newcommand{\gm}{\mbox{\scriptsize{\rm gm}}}
\newcommand{\crys}{\mbox{\scriptsize{\rm crys}}}
\newcommand{\logcrys}{{\rm logcrys}}

\renewcommand{\bar}{\overline}
\newcommand{\into}{\hookrightarrow}
\renewcommand{\implies}{\mbox{$\Rightarrow$}}
\newcommand{\veq}{\mbox{\large $\parallel$}}  %For vertical `equals'
\newcommand{\sZ}{\mbox{\scriptsize{$\Z$}}}   %For \Z (integers)
\newcommand{\sC}{\mbox{\scriptsize{$\C$}}}   %\C in subscript
\newcommand{\sQ}{\mbox{\scriptsize{$\Q$}}}   %\Q in subscript
\newcommand{\limdir}[1]{\mathop{\rm
lim}_{\buildrel\longrightarrow\over{#1}}}
\newcommand{\liminv}[1]{\mathop{\rm
lim}_{\buildrel\longleftarrow\over{#1}}}
\newcommand{\onto}{\mbox{$\to\!\!\!\!\to$}}

\newcommand{\boxtensor}{\def\boxtimesten{\Box\kern-7.59pt\raise1.2pt
\hbox{$\times$} }}                                  %For boxtensor

\newcounter{elno}                   % This to number lists

\newcommand{\cC}{\mathcal{C}}

\newcommand{\cF}{\mathcal{F}}

\newcommand{\cL}{\mathcal{L}}

\newcommand{\cO}{\mathcal{O}}
\newcommand{\cP}{\mathcal{P}}

\newcommand{\cT}{\mathcal{T}}

\newcommand{\cZ}{\mathcal{Z}}

\begin{document}

\title{On the Theory of 1-Motives}
\author{Luca Barbieri-Viale}
\address{Dipartimento di Matematica Pura e Applicata, Universit\`a degli Studi di Padova\\ Via G. Belzoni, 7\\Padova -- I-35131\\ Italy}
\email{barbieri@math.unipd.it}
\date{December 12,  2005}
\begin{abstract} This is an overview and a preview of the theory of {\it mixed motives of level $\leq$ 1}\, explaining some results, projects, ideas and indicating a bunch of problems. 
\end{abstract}
\maketitle
\hfill {\it\small Dedicated to Jacob Murre}\\

Let $k$ be an algebraically closed field of characteristic zero to start with and let $S= \Spec (k)$ denote our base scheme.
Recall that Murre \cite{MU} associates to a smooth $n$-dimensional projective variety $X$ over $S$ a Chow cohomological {\it
Picard motive} $M^1(X)$ along with the {\it Albanese motive} $M^{2n-1}(X)$. The projector $\pi_1\in CH^n(X\times X)_{\sQ}$ defining $M^1(X)$ is obtained
{\it via}\, the isogeny $\Pic^0 (X) \to \Alb (X) $ between the Picard and Albanese variety, given by the restriction to a smooth curve $C$ on $X$ since $\Alb (C) = \Pic^0 (C)$ (such a curve is obtained by  successive hyperplane sections). 
For a survey of classical Chow motives see \cite{SC} (\cf also \cite{MJM}).\\

In the case of curves $M^1 (X)$ is the Chow motive of $X$ refined from lower and higher trivial components, \ie $M^0(X)$ and
$M^{2}(X)$, such that, for smooth projective curves $X$ and $Y$ \begin{equation}\label{main}
\Hom (M^1(X),M^1(Y)) \cong \Hom (\Pic^0 (X), \Pic^0 (Y))_{\sQ}
\end{equation} by Weil (see \cite[Thm. 22 on p. 161]{WE} and also a remark of
Grothendieck and Manin \cite{MA}). 
Furthermore, the semi-simple abelian category of abelian
varieties up to isogeny is the pseudo-abelian envelope of the category of Jacobians and $\Q$-linear maps. Thus, such a theory of {\em pure} motives of smooth projective curves is known to be equivalent to the theory of abelian varieties up to
isogeny, as pointed out by Grothendieck: {\em one-dimensional (pure) motives} are abelian varieties.\\

This formula \eqref{main} suggests that we may take objects represented by $\Pic$-functors as models for larger categories of
{\em mixed} motives of any kind of curves over arbitrary base schemes $S$. However, non representability of $\Pic$ for open
schemes, forces to refine our models. Let $\bar X$ be a closure of $X$ with divisor at infinity $X_{\infty}$, \ie $X = \bar X -
X_{\infty}$. For $\bar X$ smooth we have that $\Pic (X)$ is the cokernel of the canonical map $\Div_{\infty}(\bar X)\to \Pic
(\bar X)$ associating $D\mapsto \cO(D)$, for divisors $D$ on $\bar X$ supported at infinity. Thus we may set our models
following Deligne \cite{D} and Serre \cite{SE} as
\begin{equation}\label{open} [\Div_{\infty}^0(\bar X)\to \Pic^0 (\bar X)]\end{equation} when $X$ is smooth over $S= \Spec (k)$,
by mapping algebraically equivalent to zero divisors at infinity to line bundles. Therefore, a vague defintion of our categories
of {\em 1-motives}\, $M$ can be visioned as two terms complexes (up to quasi-isomorphisms) of the following kind
$$M\df [L \to G]$$ where $L$ is {\em discrete-infinitesimal} and $G$ is {\em continuous-connected}. Moreover, we expect that a
corresponding formula \eqref{main} would be available in the larger category of mixed motives. 

\section{On Picard functors} Let $\pi : X\to S$ and consider the {\it Picard functor}\,  $T \leadsto \Pic_{X/S}(T)$ on the
category of schemes over $S$ obtained by sheafifying the functor $T\leadsto  \Pic (X\times_S T)$ with respect to the
fppf-topology (= flat topology). This means that if $\pi:X\times_ST\to T$ then
$$\Pic_{X/S}(T)\df H^0_{\rm fppf}(T,R^1\pi_*(\G_m\mid_{X\times_ST})).$$ If $\pi_*(\cO_{X})=\cO_S$ or by reducing to this
assumption, \eg if $\pi$ is proper, the Leray spectral sequence along $\pi$ and descent yields an exact sequence
$$0\to
\Pic (S)\to \Pic (X)\to \Pic_{X/S}(S)\to H^2_{\rm fppf}(S,\G_m)\to H^2_{\rm fppf}(X,\G_m).$$ Here the \'etale topology will
suffices as $H^i_{\et}(-, \G_m)\cong H^i_{\rm fppf}(-, \G_m)$ for all $i\geq 0$ by a theorem of Grothendieck (see \cite[VI.5 p.
126 \& VI.11 p. 171]{XEX}). If there is a section of $\pi$ we then have that 
$\Pic_{X/S}(S)\cong \Pic (X)/\Pic (S).$

If we set $\pi : X\to S$ proper and flat over a base, the Picard fppf-sheaf $\Pic_{X/S}$ would be possibly representable by an
algebraic space only. For a general theory we should stick to algebraic spaces not schemes (see \cite[8.3]{BLR}). However, as
far as $S= \Spec (k)$ is a field, we may just consider group schemes:  by Grothendieck and Murre (see \cite{MUR} and \cite[8.2]{BLR}) we
have that $\Pic_{X/k}$  is representable by a scheme locally of finite type over $k$. As a group scheme $\Pic^0$ usually stands
for the connected component of the identity of $\Pic$ and $\Pic^0_{X/k}$ is an abelian variety (known classically as the Picard
variety, \cf \cite[8.4]{BLR}) as soon as $X$ is also smooth and $k$ has zero characteristic. Here $\NS_X \df \pi_0 (\Pic_{X/k})$
is finitely generated. In positive characteristic, for $X$ smooth and proper over $k$ perfect, the connected component of the
identity endowed with its reduced structure $\Pic^{0,\red}_{X/k}$ is an abelian variety. More informations, \eg on the universal
line bundle $\cP$ on $X\times_S \Pic_{X/S}$, can be obtained from \cite[\S 8]{BLR}.

For example, if $X$ is a singular projective curve, in zero characteristic, and $\tilde X$ is the normalization of $X$ we then
have an extension\footnote{In the geometric case, \ie when $k$ is algebraically closed, $\Pic_{X/k}(k)\cong \Pic(X)$.}  
$$0\to V\oplus T\to \Pic^0 (X)\to \Pic^0(\tilde X)\to 0$$ where $V=\G_a^{r}$ is a vector group and $T=\G_m^{s}$ is a torus. The
additive part here is non homotopical invariant, that is, the semi-abelian quotient is homotopical invariant, \eg consider the
well known example of $X$ = projective rational cusp: its  first singular cohomology group is zero but $\Pic^0_{X/k} = \G_a$. 
For proper schemes in zero characteristic, we can describe the semi-abelian quotient of $\Pic$-functors as follows. 

\subsection{Simplicial Picard functors} Let $\pi :X_{\d}\to X$ be a smooth proper hypercovering of $X$ over $S =\Spec (k)$.
Recall that $X_{\d}$ is a simplicial scheme with smooth components obtained roughly as follows: $X_0$ is a resolution of
singularities of $X$, $X_1$ is obtained by a resolution of singularities of $X_0\times_X X_0$, etc. Such hypercoverings were
introduced by Deligne \cite{D} in characteristic zero (after Hironaka's resolution of singularities) but are also available over
a perfect field of positive characteristic (after de Jong's theory \cite{DJ}) by taking $X_0$ an alteration of $X$ (in this case
$X_0\to X$ is only generically \'etale). Actually, it is possible to refine such a construction, in characteristic zero,
obtaining a (semi)simplicial scheme $X_{\d}$ such that $\dim (X_{i})= \dim (X) -i$ so that the corresponding complex of algebraic
varieties (in the sense of \cite{BRS}) is bounded.
\subsubsection{} Denote $\bPic (X_{\d})\cong\HH^1_{\rm fppf}(X_{\d},\cO_{X_{\d}}^*)\cong \HH^1_{\et}(X_{\d},\cO_{X_{\d}}^*)$ the
group of isomorphism classes of simplicial line bundles on $X_{\d}$, \ie of invertible
$\cO_{X_{\d}}$-modules. Let $\bPic_{X_{\d}/S}$ be the associated fppf-sheaf on $S$. Over $S=\Spec (k)$ such $\bPic_{X_{\d}/S}$ 
is also representable (see \cite[A]{BSAP}). The canonical spectral sequence for the components of $X_{\d}$ yields the following
long exact sequence of fppf-sheaves:
\begin{equation}\label{semisimp} 
\begin{array}{r}
 \frac{\ker ((\pi_1)_*\G_{m,X_1}\to (\pi_2)_*\G_{m,X_2})}{\im ((\pi_0)_*\G_{m,X_0}\to (\pi_1)_*\G_{m,X_1})}\into
\bPic_{X_{\d}/S} \to \ker (\Pic_{X_0/S}\to\Pic_{X_1/S})\\[10pt]
\to \frac{\ker ((\pi_2)_*\G_{m,X_2}\to (\pi_3)_*\G_{m,X_3})}{\im ((\pi_1)_*\G_{m,X_1}\to (\pi_2)_*\G_{m,X_2})}
\end{array}
\end{equation} where $\pi_i :X_i\to S$ are the structure morphisms.  By pulling back along $\pi :X_{\d}\to X$ we have the
following natural maps $$\Pic_{X/S} \by{\pi^*} \bPic_{X_{\d}/S}\to  \ker (\Pic_{X_0/S}\to\Pic_{X_1/S}).$$ The most wonderful
property of hypercoverings is {\it cohomological descent} that is  an isomorphism $$H^*_{\et}(X,\cF)\cong
\HH^*_{\et}(X_{\d},\pi^*(\cF))$$ for any sheaf  $\cF$ on $S_{\et}$ (as well as for other usual topologies).
 In particular, for the \'etale sheaf $\mu_m\cong \Z/m$ of $m$-rooths of unity on $S=\Spec (k)$, $k=\bar k$ and $(m,{\rm char}
(k))=1$, by (simplicial) Kummer theory (see \cite[5.1.2]{BSAP}) and cohomological descent we get the following commutative
square of isomorphisms
$$\begin{array}{ccc}
\HH^1_{\et}(X_{\d},\mu_m)&\cong&\bPic (X_{\d})_{m-\tor}\\
\veq\quad & &\veq\quad \\ H^1_{\et}(X,\mu_m)&\cong&\Pic (X)_{m-\tor}
 \end{array}$$ The simplicial N\'eron-Severi group $\NS (X_{\d})\df \bPic (X_{\d})/\bPic^0 (X_{\d})$ is finitely generated,
therefore the Tate module of $\bPic (X_{\d})$ is isomorphic to that of $\bPic^0(X_{\d})$ and, by cohomological descent, to that
of $\Pic^0(X)$.  Moreover,  $\bPic^0(X_{\d})$ is the group of $k$-points of a semi-abelian variety, in which torsion points are
Zariski dense.
\subsubsection{\bf Scholium (\protect{\cite[5.1.2]{BSAP}})}\label{surj} {\em If $X$ is proper over $S=\Spec (k)$, $k=\bar k$ of
characteristic 0, and $\pi:X_{\d}\to X$ is any smooth proper hypercovering, then 
$$\pi^*:\Pic^0(X)\onto\bPic^0(X_{\d})$$ is a surjection with torsion free kernel.\\ }

As a consequence, we see that the simplicial Picard variety
$\bPic^0(X_{\d})$ is the semi-abelian quotient of the connected commutative algebraic group $\Pic^0(X)$. Moreover, if $X$ is
semi-normal, then $\pi_*(\cO_{X_{\d}}^*)=\cO_{X}^*$, and so
$\pi^*:\Pic(X)\into\bPic(X_{\d})$ is injective, by the Leray spectral sequence for the sheaf
$\cO_{X_{\d}}^*$ along $\pi$; therefore, from~\ref{surj} we get
$$\Pic^0(X)\cong\bPic^0(X_{\d})\cong \ker^0 (\Pic^0(X_0)\to\Pic^0(X_1))$$ whenever $\Pic(X)\into\Pic (X_{0})$ is also injective
(here $\ker^0$ denotes the  connected component of the identity of the kernel). Thus, if $X$ is normal $\Pic^0(X)$ is an abelian
variety which can be represented in terms of $X_0$ and $X_1$ only. If $X$ is only semi-normal a similar argument applies and
$\Pic^0(X)\cong \bPic^0(X_{\d})$ is semi-abelian.

\subsubsection{} Homotopical invariance of units and Pic, \ie $H^i (X, \G_m)\cong H^i (\Aff^1_S\times_S X, \G_m)$ for $i=0,1$ induced by the
projection $\Aff^1_S\times_S X \to X$, is easily deduced for $X$ smooth.  Let  $\Aff^1_S\times_S X_{\d}\to X_{\d}$ be the
canonical projection; considering  \eqref{semisimp} we see that 
$$\bPic (X_{\d})\cong \bPic (\Aff^1_S\times_S X_{\d})$$ since $X_{\d}$ has smooth components.  Therefore, the semi-abelian
quotient of $\Pic^0 (X)$ is always homotopical invariant. By dealing with homotopical invariant theories we just need to avoid
the additive factors, and $\bPic^0$ is the `motivic' object corresponding to $M^1$ of proper (arbitrarily singular) $S$-schemes, \ie $\Pic^+$ in the notation adopted in \cite{BSAP} (\cf  \cite{RA1} and also the commentaries below \ref{dconj}).

\subsubsection{} In positive characteritic $p>0$ the picture is more involved and a corresponding {\em Scholium}~\ref{surj} is
valid up to $p$-power torsion only. However, the semi-abelian scheme $\bPic^{0,\red}(X_{\d})$ is independent of the choices of
the hypercovering $X_{\d}$ (see \cite[A.2]{AB}) furnishing a motivic definition of 
$H^1_{\crys}$ (described in \cite{AB}, \cf \ref{crysconj} below). 

\subsection{Relative Picard functors} For a pair $(X, Y)$ consisting of a proper $k$-scheme $X$ and a closed sub-scheme $Y$ we
have a natural long exact sequence
\begin{equation}\label{longpic} H^0(X,\cO_{X}^*)\to H^0(Y,\cO_{Y}^*)\to
\Pic (X, Y)\to \Pic (X)\to \Pic (Y)
\end{equation} induced by the surjection of Zariski (or fppf) sheaves
$\G_{m,X}\to i_*\G_{m, Y}$ where $i: Y\into X$ is the inclusion; here
\[\Pic (X, Y) =\HH^1(X,\G_{m,X}\to i_*\G_{m,Y})\] is the group of isomorphism classes of  pairs $(\cL, \varphi)$ such that
$\cL$ is a line bundle on $X$ and  $\varphi :\cL\mid_{Y}\cong \cO_{Y}$ is a trivialization on $Y$ (see \cite[\S 2]{BSAP}). For
$(X, Y)$ as above the fppf-sheaf associated to the relative Picard functor
$$T\leadsto  \Pic (X\times_k T, Y\times_k T)$$ is representable by a $k$-group scheme which is locally of finite type over $k$
(\cf \cite[A]{BSAP}). If $\Pic^0(X)$ is abelian, \eg $X$ is normal, the sequence (\ref{longpic}) yields a semi-abelian group
scheme $\Pic^0(X, Y)$  (\cf \cite[2.1.2]{BSAP})  which can be represented as an extension (say of $k$-points over $k= \bar k$ of
characteristic zero)
\begin{equation}
\frac{ H^0(Y,\cO_{Y}^*)}{\im H^0(X,\cO_{X}^*)}\into \Pic^0(X, Y)\onto 
\ker^0(\Pic^0(X)\to\Pic^0(Y))
\end{equation} where $\Pic^0(X, Y)$ is the connected component of the identity of $\Pic (X,Y)$, the $k$-torus is
$\coker\left((\pi_{X})_*\G_{m,X} \to (\pi_{Y})_*\G_{m,Y}\right)$ where
$\pi_{X}: X\to\Spec k$, $\pi_{Y}:Y\to\Spec k$ are the structure morphisms and where $\ker^0$ denotes the connected component of
the identity of the kernel (the abelian quotient is further described below).

\subsubsection{} For example, assume $X$ proper (normal) and $Y=\cup Y_i$, where $Y_i$ are the (smooth) irreducible components
of a reduced normal crossing divisor $Y$. Consider the normalization $\pi :\coprod Y_i\to Y$ and observe that $\pi^*:\Pic (Y)\to
\oplus\Pic (Y_i)$ is representable by an affine morphism (see \cite[2.1.2]{BSAP}). Therefore $$\ker^0(\Pic^0(X)\to\Pic^0(Y)) =
\ker^0(\Pic^0(X)\to\oplus\Pic^0(Y_i)).$$ Moreover, for any such pair $(X, Y)$,  we have that (\cf \cite[2.2]{BSAP}) any relative
Cartier divisor  $D\in\Div (X, Y)$, \ie as divisor on $X$ such that the support $|D|\cap Y =\emptyset$,  provides $(\cO_{X}(D),
1)$ which defines an element  $[D]\in\Pic (X, Y)$ where $1$ denotes the tautological section of $\cO_{X}(D)$, trivializing it on
$X - |D|$. Here a Cartier divisor $D\in\Div (X, Y)$ is algebraically equivalent to zero relative to $Y$ if $[D]\in \Pic^0(X,
Y)$. Denote
$\Div^0_Z(X, Y)\subset \Div_Z (X, Y)$ the subgroup of relative divisors supported on a closed sub-scheme  $Z\subset X$ which are
algebraically equivalent to zero relative to $Y$. We also have a `motivic' object
$$ [\Div^0_Z(X, Y)\to \Pic^0(X, Y)]$$ which morally corresponds to $M^1(X-Z,Y)$.

\subsubsection{} Starting from an open scheme $X$ let $\bar X$ be a closure of $X$ with boundary $X_{\infty}$, \ie 
$X = \bar X - X_{\infty}$. For $Z=X_{\infty}$ and $Y=\emptyset$ from the pair $(\bar X, \emptyset)$ we get  \eqref{open} and for
$Y=X_{\infty}$ we have $[\Div^0_Z(\bar X, X_{\infty})\to \Pic^0(\bar X, X_{\infty})]$ (\cf \cite[2.2.1]{BSAP}).

\subsection{Higher Picard functors} Let $X$ be an equidimensional $k$-scheme. Let $$CH^p(X)\df \cZ^p(X)/\equiv_{\rm rat}$$ be the Chow group of codimension 
$p$-cycles modulo rational equivalence. Recall that $CH^1(X) =\Pic (X)$ if $X$ is smooth but the Chow functor $T\leadsto CH^p (X\times_k T)$ for $1<p\leq \dim (X)$ doesn't provide a representable functor even in the case when $X$ is smooth and proper over $k =\bar k$.

\subsubsection{} \label{CH}
To supply this defect several proposed generalizations have been investigated (see  \cite{MU1}, \cite{LB} and \cite{GR}).
Consider the sub-group $CH^p(X)_{\rm alg}$ of those cycles in $CH^p(X)$ which are algebraically equivalent to zero and let $NS^p(X)\df CH^p(X)/CH^p(X)_{\rm alg}$ denote the N\'eron-Severi group. Denote $CH^p(X)_{\rm ab}$ the sub-group of
$CH^p(X)_{\rm alg}$ of those cycles which are abelian equivalent to zero\footnote{Note that $CH^p(X)_{\rm alg}$ and $CH^p(X)_{\rm ab}$ are divisible groups.}, 
\ie  $CH^p(X)_{\rm ab}$ is the intersection of all kernels of regular  homomorphisms from $CH^p(X)_{\rm alg}$ to abelian varieties (see \cite{MU1} for definitions and references). The main question here is about the existence of an `algebraic representative', \ie a universal regular homomorphism from $CH^p(X)_{\rm alg}$ to an abelian variety. In modern terms, one can rephrase it (equivalently or not) by asking if the homotopy invariant sheaf with transfers  (see \cite{V} for this notion) $CH_{X/k}^p$ associated to $X$ smooth is provided with a universal map to a {\em 1-motivic sheaf} (see \cite{BK} and \cite{AY}, also \ref{Tot} below). The abelian category $Shv_1(k)$ of 1-motivic \'etale sheaves is given by those homotopy invariant sheaves with transfers $\cF$ such that there is map $$ G \by{f} \cF$$ where $G$ is {\em continuous-connected} (\eg semi-abelian); $\ker f$ and $\coker f$ are {\em discrete-infinitesimal} (\eg finitely generated). The paradigmatic example is $\cF = \Pic_{X/k}$ for $X$ a smooth $k$-variety (see \cite{BK}). Starting from $CH_{X/k}^p$ we may seek for
$$c^p: CH_{X/k}^p \to (CH_{X/k}^{p})^{(1)}$$
with $(CH_{X/k}^{p})^{(1)}\in Shv_1(k)$ universally. 
Remark that the key point is to provide a {\it finite type}\, object as such a universal  Ind-object always exists (see \cite{AY}). Namely, $CH^p(X)_{\rm alg}^{(1)}$ will be related to the `algebraic representative'.

\subsubsection{} Assume the existence of a universal regular homomorphism $\rho^p :CH^p(X)_{\rm alg}\to A^p_{X/k}(k)$ to (the
group of $k$-points) of an abelian variety $A^p_{X/k}$ defined over the base field $k$. This is given by Murre's theorem for $p=2$ (see \cite{MU1}) and it is clear for $p=1,\dim(X)$ by the theory of the Picard and Albanese varieties. We then quote the following functorial algebraic filtration $F_a^*$ on $CH^p(X)$ (\cf \cite{BA}):
\begin{itemize}
\item $F_a^0CH^p(X) = CH^p(X)$,
\item $F_a^1CH^p(X) = CH^p(X)_{\rm alg}$ 
\item $F_a^2CH^p(X) = CH^p(X)_{\rm ab}$, \ie is the kernel of the universal regular homomorphism $\rho^p$ above, 
\item  and the corresponding extension
\begin{equation}\label{NSext}
0 \to A^p_{X/k}(k) \to CH^p(X)/F^2_a \to NS^p(X)\to 0.\\
\end{equation}
\end{itemize}  

Remark that Bloch, Beilinson and Murre (see \cite{JM}) conjectured the existence of a  finite filtration $F_m^*$ on $CH^p(X)_{\sQ}$ (with rational coefficients)
such that $F_m^1CH^p(X)$ is given by $CH^p(X)_{\rm hom}$, \ie by the sub-group of those codimension $p$ cycles which are
homologically equivalent to zero for some Weil cohomology theory, $F_m^*CH^p(X)$ should be functorial and compatible with the
intersection pairing. The motivic filtration $F_m^*$ will be inducing the algebraic (or 1-motivic) filtration $F_a^*$ somehow, \eg $F_a^* = F_m^*\cap CH^p(X)_{\rm alg}$
for $*>0$.

Remark that we may even push further this picture by seeking for the 1-motivic algebraically defined extension of codimension $p$ cycles modulo numerical equivalence by $A^p_{X/k}(k)$  (which pulls back to \eqref{NSext}, see also \ref{univab} below). 

\section{On 1-motives} A {\it free}\,  1-motive over $S$ (here $S$ is any base scheme) in Deligne's definition (a {\it 1-motif
lisse}\, \cf \cite[\S 10]{D}) is a complex $M \df [L \by{u} G]$ of $S$-group schemes where $G$ is semi-abelian, \ie it is an
extension of an abelian scheme $A$ by a torus $T$ over $S$, the group scheme $L$ is, locally for the
\'etale topology on $S$, isomorphic to a finitely-generated {\it free}\, abelian constant group, and $u: L\to G$ is an
$S$-homomorphism. A 1-motive $M$ can be represented in a diagram
$$ \begin{array}{ccc}
 & L & \\ & \downarrow & \\
 0\to \ T& \to\ G\ \to & A\ \to 0
\end{array}$$

An {\em effective}\, morphism of 1-motives is a morphism of the corresponding complexes of group schemes  (and actually of the
corresponding diagrams). Any such a complex can be regarded as a complex of fppf-sheaves. Following the existing literature, $L$
is regarded in degree $-1$ and $G$ in degree $0$ (however, for some purposes, \eg in order to match the conventions in Voevodsky
triangulated categories, it is convenient to shift $L$ in degree $0$ and $G$ in degree $1$, \cf \cite{BK}). We let $\M^{\fr}$ denote the category of Deligne 1-motives (\cf \ref{torsion} below).

\subsection{Generalities} It is easy to see that $\M^{\fr}$ has kernels and cokernels but images and coimages,  in general,
don't coincide. For kernels, if 
$\ker^c (\phi) =[\ker (f) \by{u}  \ker (g)]$ is the
kernel of $\phi=(f,g):M\to M'$ as a map of complexes then 
$\ker (\phi)= [\ker^0(f)
\by{u}  \ker^0(g)]$ is the pull-back of
$\ker^0(g)$ along $u$, where $ \ker^0(g)$ is the connected component 
of the identity of the kernel  of $g : G\to G'$ and $\ker^0(f) \subseteq\ker(f)$.

Similarly, for cokernels, if 
$\coker^c (\phi) =[\coker (f)\by{
\bar u'}\coker(g)]$ is the cokernel as complexes and $\cT$ is the torsion subgroup of $\coker (f)$, as group schemes, then $$\coker
(\phi) = [\coker (f)/\cT \to \coker(g)/\bar u' (\cT) ]$$ is a Deligne's 1-motive which is clearly a cokernel of $\phi$.

Associated to any 1-motive $M$ there is a canonical extension (as two terms complexes) 
\begin{equation}\label{motexseq} 0\to [0\to G] \to M\to [L\to 0]\to 0
\end{equation}
\subsubsection{} Actually, a 1-motive $M$ is canonically equipped with an increasing {\it weight filtration}\,  by sub-1-motives
as follows:
$$W_i(M) =\left \{\begin {array}{rr} M & i\geq 0\\[4pt] [0\to G]& i= -1\\[4pt] [0\to T] & i= -2\\[4pt] 0 & i \leq -3 \end{array}
\right. $$ In particular we have $\gr_{-1}^W(M) = [0\to A]$ and $\gr_0^W(M)=[L\to 0]$.\\
\subsubsection{} For $S = \Spec (k)$ a 1-motive $M =[L\by{u} G]$ over $k$ (a perfect field)  is equivalent to the given semi-abelian $k$-scheme $G$, a finitely generated free abelian  group $\bar
L$ which underlies a $\Gal(\bar k/k)$-module, a 1-motive $[\bar L\by{\bar u} G_{\bar k}]$ over $\bar k$, such that
$\bar u$ is $\Gal(\bar k/k)$-equivariant, for the given module structure on $\bar L$, and the natural semi-linear action on
$G_{\bar k}=G\times_k\bar k$. In fact, the morphism $u$ is determined uniquely by base change to $\bar k$, \ie by the morphism 
$u_{\bar k}:L_{\bar k}\to G_{\bar k}$, which is $\Gal(\bar k/k)$-equivariant. 
\subsubsection{} It is easy to see that there are no non-trivial quasi-isomorphisms between Deligne 1-motives.  Actually, there
is a canonical functor $\iota : \M^{\fr} \to D^b(S_{\rm fppf})$ which is a full embedding into  the derived category of bounded
complexes of sheaves for the fppf-topology on $S$.

\subsubsection{\bf Scholium (\protect{\cite[Prop.2.3.1]{RAY}})}\label{freein}  {\em Let $M$ and $M'$ be free 1-motives. Then
$$\Hom_{\M^{\fr}} (M, M')\cong \Hom_{D^b(S_{\rm fppf})} (\iota (M), \iota (M')).$$}
\begin{proof} The naive filtration of $M =[L\to G]$ and $M'=[L\to G']$ yields a spectral sequence 
$$E_1^{p, q} = \bigoplus_{-i+j =p} \Ext^q ({}^iM, {}^jM')\implies
\EExt^{p+q}(M, M')$$ yielding complexes $E_1^{\d , q}$
$$ \Ext^q (G, L')\to \Ext^q (G, G')\oplus \Ext^q (L, L')\to 
\Ext^q (L, G'),$$ where the left-most non-zero term is in degree -1. We see that $\Ext^0(G, L')\allowbreak =\Hom (G, L')=0$ since
$G$ is connected and $\Ext(G, L')=0$ since $L'$ is free. Thus
$E^{0,0}_2 = \Hom_{\M^{\fr}} (M, M')$ is 
$\EExt^0(M, M')= \Hom_{D^b(S_{\rm fppf})} (\iota (M), \iota (M')).$
\end{proof}

\subsection{Hodge realization} \label{Hodge} The {\it Hodge realization} $T_{Hodge}(M)$ of a 1-motive $M$ over $S=\Spec (\C)$
(see \cite[10.1.3]{D}) is $(T_{\sZ}(M), W_*, F^*)$ where $T_{\sZ}(M)$ is  the lattice given by the pull-back of  $u : L\to G$
along $\exp : \Lie (G) \to G$, $W_*$ is the integrally defined weight filtration\footnote{Note that $H_1(G)$ is the kernel of $\exp :  \Lie (G) \to G$.}
$$W_iT(M) \df\left \{\begin {array}{rr} T_{\sZ}(M) & i\geq 0\\ H_1(G) & i= -1\\ H_1(T) & i= -2\\ 0 & i \leq -3
\end{array} \right.$$ and $F^*$ is the Hodge filtration defined by
$F^0(T_{\sZ}(M)\otimes\C)\df \ker (T_{\sZ}(M)\otimes\C\to \Lie(G)).$ Then we see that $T_{Hodge}(M)$ is a mixed Hodge structure and we have
$\gr_{-1}^WT_{Hodge}(M)\cong H_1(A,\Z)$ as pure polarizable Hodge structures of weight
$-1$. 
\subsubsection{} The functor
$$T_{Hodge}: \M^{\fr}(\C)\longby{\simeq} {\rm MHS}_1^{\fr}$$ is an equivalence between the category of 1-motives over $\C$ and
the category of torsion free $\Z$-mixed Hodge structures of type $$\{(0,0), (0,-1), (-1,0), (-1,-1)\}$$ such that $\gr_{-1}^W$
is polarizable. Deligne (\cf \cite[\S 10.1.3]{D}) observed that such a $H\in {\rm MHS}_1^{\fr}$
is equivalent to a 1-motive over the complex numbers. In fact, for $H\in {\rm MHS}_1^{\fr}$ the canonical extension of mixed
Hodge structures
\begin{equation}\label{hodgeexseq} 0\to W_{-1}(H)\to H\to \gr_{0}^W(H)\to 0
\end{equation} yields an {\em extension class map}\, (\cf \cite{CA})
$$e_H : \Hom_{\rm MHS} (\Z, \gr_{0}^W(H))\to \Ext_{\rm MHS}(\Z, W_{-1}(H))$$ which provides a 1-motive with lattice $L \df
\gr_{0}^W(H_{\sZ})$ mapping to the semi-abelian variety with complex points $G(\C)\df\Ext_{\rm MHS} (\Z, W_{-1}(H))$. Summarizing up, any
1-motive $M$ over $\C$ has a {\em covariant}\,  Hodge realization
$$M\leadsto T_{Hodge}(M)$$ and the exact sequence (\ref{motexseq}) gives rise to the exact sequence \eqref{hodgeexseq} of Hodge realizations.
\subsubsection{} We have that $T_{Hodge}([0\to \G_m])=\Z (1)$ is the Hodge structure (pure of weight $-2$ and purely of type
$(-1,-1)$) provided by the complex exponential  $\exp : \C\to \C^*$, \ie here $T_{\sZ}([0\to \G_m])$ is the free $\Z$-module on
$2\pi\sqrt{-1}$. Recall that for $H\in {\rm MHS}_1^{\fr}$ we get $H^\vee\df\ihom (H, \Z (1))\in {\rm MHS}_1^{\fr}$ where $\ihom$
is the internal Hom in ${\rm MHS}$ (see \cite{D}). We have that $\Z^{\vee} = \Z (1)$. Moreover 
$$(\ \ )^\vee \df \ihom (\ \ , \Z (1)): ({\rm MHS}_1^{\fr})^{op}\by{\simeq} {\rm MHS}_1^{\fr}$$ is an anti-equivalence providing
${\rm MHS}_1^{\fr}$ of a natural involution. We may set a {\em contravariant}\, Hodge realization given by
$$M\leadsto T^{Hodge}(M)\df T_{Hodge}(M)^{\vee}$$ and an induced involution on $\M^{\fr}(\C)$ defined by the formula $$M^{\vee}
\df T_{Hodge}^{-1}\circ T^{Hodge}(M).$$ Actually, such an involution can be made algebraic (see \ref{Cartier} below) and is
known as Cartier duality for 1-motives.
\subsubsection{} Remark that ${\rm MHS}_1^{\fr}\subset {\rm MHS}_1$ where we just drop the assumption that the underlying
$\Z$-module is torsion free and we have that the category ${\rm MHS}_1$ is a thick abelian sub-category of  (graded polarizable)
mixed Hodge structures. In \cite[\S 1]{BRS} an algebraic description of ${\rm MHS}_1$   is given (see \ref{torsion} below). For
$H\in {\rm MHS}$ let $H_{(1)}$ denote the maximal sub-structure of the considered type (= largest 1-motivic sub-structure, for
short)  and let $H^{(1)}$ be the largest 1-motivic quotient. For $H'\in {\rm MHS}_1$ we clearly have
 $$\Hom_{\rm MHS} (H',H) = \Hom_{{\rm MHS}_1} (H', H_{(1)}) $$ and
 $$\Hom_{\rm MHS} (H,H') = \Hom_{{\rm MHS}_1}  (H^{(1)}, H').$$ In other words the embedding  ${\rm MHS}_1 \subset  {\rm MHS}$
has right and left adjoints given by the functors $H\mapsto H_{(1)}$ and $H\mapsto H^{(1)}$ respectively. Moreover, it is quite
well known to the experts that $\Ext^1_{{\rm MHS}_1}$ is right exact and the higher extension groups $\Ext^i_{{\rm MHS}_1}$ ($i>1$) vanish since similar
assertions hold in ${\rm MHS}$ (by Carlson \cite{CA}) and the objects of  ${\rm MHS}_1$ are stable by extensions in ${\rm MHS}$.
As a consequence, the derived category $D^b({\rm MHS}_1)$ is a full sub-category of $D^b({\rm MHS})$.
  
\subsection{Flat, $\ell$-adic and \'etale realizations} \label{flat} 
 Let $M=[L\by{u} G]$ be a 1-motive over $S$ which we consider as a complex of fppf-sheaves over $S$ with $L$ in degree $-1$ and
$G$ in degree $0$. 
 Consider the cone $M/m$ of the multiplication by $m$ on $M$. The exact sequence (\ref{motexseq}) of 1-motives yields a short
exact sequence of cohomology sheaves
\begin{equation}\label{exactseqm}0\to H^{-1}(G/m)\to H^{-1}(M/m)\to H^{-1}(L[1]/m)\to 0\end{equation} as soon as $L$ is
torsion-free, \ie $H^{-2}(L[1]/m)=\ker (L\by{m}L)$ vanishes, since multiplication by $m$ on $G$ connected is an epimorphism of
fppf-sheaves, \ie $H^{0}(G/m)=\coker (G\by{m}G)$ vanishes. Here $H^{-1}(G/m)=$ $m$-torsion of $G$ and  $H^{-1}(L[1]/m)= L/m$
whence the sequence above is given by finite group schemes. The {\em flat realization}\, $$T_{\sZ/m}(M)\df H^{-1}(M/m)$$ is a
finite group scheme, flat over $S$, which is \'etale if $S$ is defined over $\Z [\frac{1}{m}]$. By taking the Cartier dual we
also obtain a contravariant flat realization
$T^{\sZ/m}(M)\df \shom_{\rm fppf} (H^{-1}(M/m), \G_m)$.
\subsubsection{} If $\ell$ is a prime number then the $\ell$-{\it adic realization}\, $T_{\ell}(M)$ is the inverse limit over
$\nu$ of $T_{\sZ/\ell^{\nu}}(M)$.  We have $T_{\ell}([0\to \G_m])= \Z_{\ell}(1)$ by the Kummer sequence. The $\ell$-adic
realization of an abelian scheme $A$ is the $\ell$-adic Tate module of $A$. In characteristic zero then
$$\hat{T}(M)\df\liminv{m}T_{\sZ/m}(M)=\prod_{\ell} T_{\ell}(M)$$ is called the {\em \'etale realization} of $M$. For $S=\Spec
(k)$, $\hat{T}(M_{\bar k})$, along with a natural action of $\Gal(\bar k/k)$, is a (filtered) Galois module which is a free
$\hat{\Z}$-module of finite rank. Over $S=\Spec (k)$ and $k=\bar k$ we just have $$T_{\sZ/m}(M)(k) = \frac{\{(x,g)\in L\times
G(k)\mid u(x)=-mg\}}{\{(mx,-u(x)) \mid x \in L\}}.$$
\subsubsection{} If $k=\C$ we then have a comparison isomorphism $\hat{T}(M) \cong T_{\sZ}(M)\otimes\hat{\Z}$ where $T_{\sZ}(M)$
is the $\Z$-module underlying to $T_{Hodge}(M)$ (\cf \cite[\S 1.3]{BSAP}). 

\subsection{Crystalline realization}  Let $S_0$ be a scheme and $p$ a prime number such that $p$ is locally nilpotent on $S_0$.
Now let $S_0\into S_n$ be a thickening defined by an ideal with
nilpotent divided powers.  Actually, over $S_0=\Spec (k)$ a perfect field of characteristic~$p>0$ and\/~$\W(k)$ the Witt vectors
of\/~$k$ (with the standard divided power structure\footnote{Note that for $p=2$ the standard  divided power structure of
$\W_n(k)$ is {\it
not} nilpotent.} on its maximal ideal) a thickening $S_n = \Spec (\W_{n+1}(k))$ is given by the affine scheme defined by the
truncated Witt vectors of length $n+1$ (or  equivalently by $\W (k)/p^{n+1}$).  Suppose that $M_0\df [L_0\to G_0]$ is a 1-motive
defined over $S_0$. Consider
$$M_0[p^\infty]\df \limdir{\nu} T_{\sZ/p^{\nu}}(M_0)$$ the direct limit being taken, in terms of the explicit formula above, for
$\mu\geq \nu$, by sending the class of a point~$(x,g)$ in~$L_0 \times G_0(k)$
to the class of~$(p^{\mu-\nu}x,g)$. Such $M_0[p^\infty] $ is a $p$-divisible (or
Barsotti-Tate) group and the sequence \eqref{exactseqm} yields the exact sequence
\begin{equation}\label{exactseqinfty}
0\to G_0[p^\infty] \to M_0[p^\infty]
\to L_0[p^\infty]\to 0
\end{equation}
where~$L_0[p^\infty]\df L_0\otimes\Q_p/\Z_p$. For $M_0\df [0\to A_0]$ an abelian scheme we get back the Barsotti-Tate group of $A_0$.

\subsubsection{} Let $\D$ be the {\it contravariant Dieudonn\'e functor} from the category of $p$-divisible groups over $S_0
=\Spec (k)$ to the category of $D_k$-modules, for the Dieudonn\'e ring $D_k:=\W (k)[F,V]/(FV=VF=p)$. This $\D$ is defined as the
module of homomorphisms
from the $p$-divisible group to the group of Witt covectors over $k$  and provides an anti-equivalence from the category of
$p$-divisible groups over $k$
to the category of $D_k$-modules which are finitely generated and free as
$\W (k)$-modules (see \cite{FO}).

For any such a thickening $S_0\into S_n$ the functor $\D$ can be further extended to define a {\it crystal}\,  on the nilpotent
crystalline site on $S_0$ that is (equivalently given by) the Lie algebra of the associated universal $\G_a$-extension of the
dual $p$-divisible group, by lifting it to $S_n$ (\cf \cite{MM}, \cite{AB}). Therefore, by taking $\D (M_0[p^\infty] )$ we
further obtain a filtered $F$-{\it crystal}\, on the crystalline site of\/~$S_0$, associated to the Barsotti-Tate
group~$M_0[p^\infty]$. Recall that (see \cite{AB}) the category of filtered $F$-$\W(k)$-modules consists of finitely generated
$\W(k)$-modules endowed with an increasing filtration and a $\sigma$-linear\footnote{Here $\sigma$ is the Frobenius on $\W(k)$.}
operator, the Frobenius $F$, respecting the filtration. Filtered $F$-crystals are the objects whose underlying $\W(k)$-modules
are free and there exists a $\sigma^{-1}$-linear operator, the Verschiebung $V$, such that $V \circ F=F \circ V=p$. 
\subsubsection{} The {\em crystalline realizations}\, of $M_0$ over $S_0=\Spec(k)$ are the following filtered $F$-crystals (see
\cite[\S 1.3]{AB} where are also called Barsotti-Tate  crystals of the 1-motive $M_0$ and \cf \cite{FJ} and \cite[4.7]{KT}). The contravariant one is
$$T^{\crys}(M_0)\df
\liminv{n}\D (M_0[p^\infty])(S_0\into S_n) $$and the covariant is
$$T_{\crys}\left(M_0\right)\df
\liminv{n}\D\left(M_0[p^\infty]^\vee\right)(S_0\into S_n)$$
where $M_0[p^\infty]^\vee$ is the Cartier dual. It follows from
\eqref{exactseqinfty}  that $T_{\crys}(M_0)$ admits Frobenius and Verschiebung operators and a filtration (respected by
Frobenius and Verschiebung). 
\subsubsection{} We get $T_{\crys}([0\to \G_m] )=\W(k)(1)$ which is the filtered $F$-crystal $\W(k)$, with filtration
$W_i=\W(k)$ if $i \geq -2$ and $W_i=0$ for $i<-2$ and  with the $\sigma$-linear operator $F$ given by $1 \mapsto 1$ and
the $\sigma^{-1}$-linear operator $V$ defined by $1\mapsto p$.\\

\subsection{De Rham realization}\label{DeRham} The {\em De Rham realization}\, of a 1-motive $M=[L\by{u} G]$ over a suitable
base scheme $S$ is obtained {\it via}\, Grothendieck's idea of universal $\G_a$-extensions (\cf \cite[\S 4]{MM},
\cite[10.1.7]{D} and \cite{AB}). Consider $\G_a$ as a complex of $S$-group schemes concentrated in degree
$0$. If $G$ is any $S$-group scheme such that $\shom (G, \G_a)=0$
and $\sext (G, \G_a)$ is a locally free $\cO_S$-module of finite rank, the universal
$\G_a$-extension is an extension of $G$ by the (additive dual) vector group $\sext (G,\G_a)^{\vee}$ (see \cite{MM}).

\subsubsection{} Now for any 1-motive $M =[L\by{u}G]$ over $S$, we have $\shom (M,\G_a) =0$, and by the extension
\eqref{motexseq} $\sext (M,\G_a)$ is of finite rank. Thus we obtain a {\it universal $\G_a$-extension} $M^{\natural}$, in
Deligne's notation \cite[10.1.7]{D}, where $M^{\natural}=[L\by{u^{\natural}} G^{\natural}]$ is a complex of
$S$-group schemes\footnote{Note that $G^{\natural}$ is not the universal $\G_a$-extension of $G$ unless $L = 0$.} which is an
extension of $M$ by $\sext (M,\G_a)^{\vee}$ considered as a complex in degree zero. Here we have an extension of $S$-group
schemes
$$0\to \sext (M,\G_a)^{\vee}\to G^{\natural}\to G\to 0$$ such that $G^{\natural}$ is the push-out of the universal
$\G_a$-extension of the semi-abelian scheme $G$ along the inclusion of $\sext (G,\G_a)^{\vee}$ into $\sext (M,\G_a)^{\vee}$. The
canonical map
$u^{\natural}:L\to G^{\natural}$ such that the composition
$$L\by{u^{\natural}} G^{\natural}\to \shom (L,\G_a)^{\vee}$$ is the natural evaluation map. The {\it De Rham realization}\, of
$M$ is then defined as
$$T_{DR}(M)\df\Lie G^{\natural},$$ with the {\it Hodge-De Rham filtration}\,  given by
$$F^0T_{DR}(M)\df \ker (\Lie G^{\natural}\to \Lie G)\cong
\sext(M,\G_a)^{\vee}$$
\subsubsection{} Over a base scheme on which\/~$p$ is locally nilpotent there is a canonical and functorial isomorphism (see
\cite[Prop. 1.2.8]{AB}) 
$$(M[p^\infty])^{\natural}\times_{M[p^\infty]} G[p^\infty]
\longby{\simeq} G^{\natural} \times_{G} G[p^\infty]$$ where $(M[p^\infty])^{\natural}$ also denotes the universal
$\G_a$-extension of a  Barsotti-Tate group. In particular, we have a natural isomorphism of Lie algebras
\begin{equation} \label{crysderham}
\Lie (M[p^\infty])^{\natural} \longby{\simeq}
\Lie G^{\natural}
\end{equation}
\subsubsection{} For $S_0$ a scheme such that $p$ is locally nilpotent and
$M_0 =[L_0\to G_0]$ a 1-motive over $S_0$, let~$S_0\into S$ be a locally
nilpotent pd thickening of~$S_0$.  Let $M$
and~$M'$ be two 1-motives over~$S$ lifting~$M_0$. We have proven (see \cite[\S 3]{AB}) that there is a canonical isomorphism
$M^{\natural}\cong (M')^{\natural}$ showing that the universal $\G_a$-extension is crystalline. Define the crystal of  (2-terms
complexes of) group schemes $M_0^{\natural}$  on the nilpotent crystalline
site of\/~$S_0$ as follows
$$M_0^{\natural}(S_0 \into S) \df M^{\natural}$$ which we called the {\it universal extension crystal}\, of a 1-motive (see
\cite[\S 3]{AB}).

Applying it to $M_0$ defined over $S_0= \Spec (k)$ a perfect field and $S_n = \Spec (\W_{n+1}(k))$ we see that the De Rham
realization is a crystal indeed. Actually (see \cite[\S 4]{AB} for details) the formula \eqref{crysderham} yields:
\subsubsection{\bf Scholium (\protect{\cite[Thm. A$^\prime$]{AB}})}\label{comparison} {\em There is a comparison isomorphism of
$F$-crystals $$T_{crys}(M_0) = T_{DR}(M)$$ for any (formal) lifting $M$ over $\W (k)$ of $M_0$ over $k$.}
\subsubsection{} If $k=\C$ then the De Rham realization is also compatible with the Hodge realization; we have
$$T_{DR}(M)=T_{Hodge}(M)\otimes\C$$ as bifiltered $\C$-vector spaces, \ie we have that $H_1(G^{\natural},\Z) = H_1(G,\Z)$ thus
$T_{\sC}(M)\df T_{\sZ}(M)\otimes\C\cong \Lie G^{\natural}$ and
$M^{\natural} = [L \to T_{\sC}(M)/H_1(G,\Z)]$, see \cite[\S 10.1.8]{D}.

\subsection{Paradigma} \label{para} Let $X$ be a (smooth) projective variety over $k=\bar k$. Let 
$\Pic_{X/k}$ be the Picard scheme and $\Pic^{0, \red}_{X/k}$ the connected component of the identity endowed with its reduced
structure. Recall that $\NS_X \df \pi_0 (\Pic_{X/k})$ is finitely generated and $\Pic^{0, \red}_{X/k}$ is divisible. Recall that
we always have
$$H^1_{\rm fppf}(X,\mu_n) = \Pic (X)_{n-\tor}$$ Therefore
$$T_{\ell} ([0\to \Pic^{0, \red}_{X/k}])= H^1_{\rm fppf}(X,\Z_{\ell}(1))$$ If $\ell \neq {\rm char} (k)$ then the \'etale
topology will be enough.
\subsubsection{} \label{natural} Let $\Pic^{\natural}(X)$ be the group of isomorphism classes of {\it pairs} $(\cL,\nabla)$ where $\cL$ is a
{\it line bundle}\, on $X$ and $\nabla$ is an {\it integrable connection}\, on $\cL$. In characteristic zero then there is the
following extension
$$0\to H^0(X,\Omega^1_X) \to \Pic^{\natural , 0}(X)\to \Pic^0(X)\to 0$$ where $\Pic^{\natural, 0}$ is the the subgroup of those
pairs
$(\cL ,\nabla)$ such that $\cL\in\Pic^0$. The above extension is the group of $k$-points of the universal $\G_a$-extension of
the abelian variety
$\Pic^0_{X/k}$, $\Lie \Pic^0(X) = H^1(X, \cO_X)$ and $$\Lie \Pic^{\natural, 0}(X) = H^1_{DR}(X/k)$$ as $k$-vector spaces (as
soon as the De Rham spectral sequence degenerates). Moreover, for $k =\C$, the exponential sequence grants 
$$T_{Hodge}([0\to \Pic^0_{X/\sC}]) = H^1(X, \Z(1)).$$
\subsubsection{} In general, for an abelian $S$-scheme $A$ (in any characteristics \cf \cite[\S 4]{MM}) we have
$(A^{\vee})^{\natural}= \Pic^{\natural , 0}_{A/S}$, so that the dual of $A$ has De Rham realization
$$T_{DR}([0\to A^{\vee}]) = H^1_{DR}(A/S)(1)$$ where the twist (1) indicates that the indexing of the Hodge-De Rham filtration
is  shifted by 1 (\cf \cite[\S 2.6.3]{BSAP}). 

However, for $X$ (smooth and proper) over a perfect field $k$ of characteristic $p>0$,  the $k$-vector space $H^1_{\rm DR}(X/k)$ cannot be recovered
from the Picard scheme (as remarked by  Oda \cite{ODA}). The subspace obtained  {\it via}\, the Picard scheme is closely related
to crystalline cohomology (see \cite[\S 5]{ODA}). 
\subsubsection{} Let $X$ be smooth and proper over a perfect field $k$ of characteristic $p>0$. Let $\Pic^{\crys, 0}_{X/S_n}$ be the
sheaf on the fppf site on $S_n=\Spec (\W_{n+1}(k))$ given by the functor
associating to $T$ the group of isomorphism classes of crystals of
invertible  $\cO^{\crys}_{X\times_{S_n} T/T}$-modules (which are
algebraically equivalent to~$0$ when restricted to the Zariski site). Such $\Pic^{\crys}$ is the natural substitute of the
previous functor $\Pic^{\natural}$ and we can think $H^1_{\crys}$ as $\Lie \Pic^{\crys}$ (see \cite{AB} for details, \cf
\cite{BC}). In fact, the $S_n[\varepsilon]$-points of $\Pic^{\crys}$ reducing to the identity modulo $\varepsilon$ are the
infinitesimal deformations of $\cO^{\crys}_{X/S_n}$.   For~$A_0$, an abelian variety over~$S_0=\Spec (k)$, and an abelian
scheme~$A_n$ over~$S_n$ lifting~$A_0$, the category of  crystals of invertible~$\cO^{\crys}_{A_0\times_{S_n}
S_n[\varepsilon]/S_n[\varepsilon]}$-modules over the nilpotent
crystalline site of $A_0\times_{S_0} S_0[\varepsilon]$ relative to
$S_n[\varepsilon]$ is equivalent to the category of line bundles
over~$A_n[\varepsilon]$ with integrable connection. Hence, we have an
isomorphism of sheaves over the fppf site of $S_n$
$$(A_n^\vee)^{\natural}\cong \Pic^{\natural, 0}_{A_n/S_n} \cong
 \Pic^{\crys, 0}_{A_0/S_n}$$ and passing to $\Lie$ we get a natural isomorphism
of~$\cO_{S_n}$-modules $$T^{\crys}(A_0)\otimes
\cO_{S_n}\cong \Lie (A_n^\vee)^{\natural} \cong
\Lie \Pic^{\crys, 0}_{A_0/S_n}\cong H^1_{\crys}(A_0/S_n).$$
\subsubsection{} By applying the previous arguments to the Albanese variety $\Alb (X) = (\Pic^{0, \red}_{X/k})^{\vee}= A_0$  we
see that $\Lie \Pic^{\crys, 0} (\Alb (X))$ can  be identified to the Lie algebra of the universal extension of a
(formal) lifting of $\Pic^{0, \red}_{X/k}$ to the Witt vectors. The Albanese mapping is further inducing a canonical
isomorphism\footnote{Note that $H^0(\Alb (X),\Omega^1_X)\neq H^0(X,\Omega^1_X)$ in general, in positive characteristics.} (\cf \cite[II.3.11.2]{IL},
\cite{BC} and \cite{AB}) 
$$\Lie \Pic^{\crys, 0}(\Alb (X))\longby{\simeq} \Lie \Pic^{\crys, 0} (X).$$

In conclusion, we have
$$T_{\crys}([0\to \Pic^{0, \red}_{X/k}])\cong H^1_{\crys}(X/\W (k))$$
for $X$ a smooth proper $k$-scheme.

\subsection{Cartier duality} \label{Cartier} For $H=T_{Hodge}(M)$, $H^{\vee}=\ihom (H,\Z (1))$ is an implicit definition  (see
\ref{Hodge}) of the dual $M^{\vee}$ of a 1-motive $M$ over $\C$. In general, Deligne \cite[\S 10.2.11--13]{D} provided an
extension of Cartier duality to (free) 1-motives showing that is compatible with such Hodge theoretic involution. The main deal
here is the yoga of Grothendieck biextensions (see \cite{MUB}, \cite[VII 2.1)]{SGA7} and \cite[\S 10.2.1]{D}). 
\subsubsection{} A Grothendieck (commutative) biextension $P$ of $G_1$ and $G_2$ by $H$ is an $H$-torsor on $G_1\times G_2$
along with a structure of compatible isomorphisms of torsors $P_{g_1,g_2}P_{g_1',g_2} \cong P_{g_1g_1',g_2}$ and
$P_{g_1,g_2}P_{g_1,g_2'} \cong  P_{g_1,g_2g_2'}$ (including associativity and commutativity) for all points $g_1,g_1'$ of $G_1$
and $g_2, g_2'$ of $G_2$. Recall that an isomorphism class of a Grothendieck biextension (as commutative groups in a
Grothendieck topos) can be essentially translated by the formula (see \cite[VII 3.6.5]{SGA7})
$$ \Biext (G_1,G_2;H) = \EExt  (G_1\oo^L G_2,H).$$ Here we further have $\EExt  (G_1\oo^L G_2,H)= \EExt (G_1,\RHom (G_2,H))$ and
the canonical spectral sequence
$$E_2^{p,q} =  \Ext^p (G_1,\sext^q (G_2,H))\ \implies\  \EExt^{p+q} (G_1,\RHom (G_2,H))$$ yields an exact sequence of low degree
terms
$$\begin{array}{r} 0\to  \Ext (G_1,\shom (G_2,H))\to \Biext (G_1,G_2;H) \to 
\Hom (G_1, \sext (G_2,H))\\[4pt]
\to  \Ext^2 (G_1, \shom (G_2,H))\end{array}$$ If $\shom (G_2,H) =0$ then $\partial:\Biext (G_1,G_2;H) \cong \Hom (G_1, \sext
(G_2,H))$. In particular, for $H=\G_m$ and $G_2=A$ an abelian scheme, since $A^{\vee} = \sext (A,\G_m)$ for abelian schemes,
this isomorphism $\partial$ reduces to the more classical isomorphism (\cf \cite[4.1.3]{BK} and \ref{avatar} below).
\begin{equation}\label{repab}
\Hom (-,A^{\vee})\longby{\simeq} \Biext (-,A;\G_m) \end{equation} given by $f\mapsto (f\times 1)^*\cP_A^t$ pulling back the
(transposed) Poincar\'e $\G_m$-bixetension $\cP_A$ of $A$ and $A^{\vee}$, \ie the functor
$\Biext (-,A;\G_m)$ is representable by the dual abelian scheme.

If $G_1$ and $G_2$  are semi-abelian schemes we further have $$\Biext (A_1,A_2;\G_m)\cong\Biext (G_1,G_2;\G_m)$$ by pullback
from the abelian quotients $A_1$ and $A_2$ (see \cite[VIII 3.5-6]{SGA7}). Actually, we can regard biextensions of smooth
connected group schemes (over a perfect base field) $G_1$ and $G_2$ by $\G_m$ as invertible sheaves on $G_1\times G_2$
birigified with respect to the identity sections (see  \cite[VIII 4.3]{SGA7}).

\subsubsection{} Now let $M_i = [L_i\by{u_i} G_i]$ for $i=1, 2$ be two 2-terms complexes of sheaves.  A {\it biextension}\, $(P,
\tau, \sigma)$ of $M_1$ and $M_2$ by an abelian sheaf $H$ is given by {\it (i)} a Grothendieck biextension $P$ of $G_1$ and
$G_2$ by $H$ and a pair of compatible trivializations, \ie {\it (ii)} a biadditive section $\tau$ of the biextension $(1\times
u_2)^*(P)$ over $G_1\times L_2$, and {\it (iii)} a biadditive section $\sigma$ of the biextension $(u_1\times 1)^*(P)$ over
$L_1\times G_2$, such that {\it (iv)} the two induced sections $\tau_{\mid_{L_1\times L_2}}= \sigma_{\mid_{L_1\times L_2}}$
coincide.

Let $ \Biext (M_1,M_2;H)$ denote the group of isomorphism  classes of  biextensions. We still have the following fundamental
formula (see \cite[\S 10.2.1]{D})
$$\Biext (M_1,M_2;H) = \EExt (M_1\oo^L M_2,H) $$ here $\EExt (M_1\oo^L M_2,H) = \EExt (M_1,\RHom (M_2,H))$ where $M_i$ is
considered a complex of sheaves concentrated in degree $-1$ and $0$.

\subsubsection{} Let $M=[L\by{u}G]$ be a 1-motive where $G$ is an extension of an abelian scheme $A$ by a torus $T$. The main
goal is that the functor on 1-motives
$$N\mapsto \Biext (N,M;\G_m)$$ is representable, \ie there is a {\em Cartier dual} $M^{\vee}= [T^{\vee}\by{u^\vee}G^u]$ such that
\begin{equation}\label{repcar}
\Hom (N, M^{\vee})\longby{\simeq}\Biext (N,M;\G_m)
\end{equation} is given by pulling back the Poincar\'e biextension generalizing \eqref{repab}. More precisely, it is given by
$\varphi\mapsto (\varphi\times 1)^*\cP_M^t$ where the Poincar\'e $\G_m$-biextension $\cP_M$ is simply obtained from that of $A$
and $A^{\vee}$ by further pullback to $G$ and $G^u$ according to the above (and below) description. See \cite[10.2.11]{D} and
\cite[1.5]{BSAP} for the construction of $M^{\vee}$ and \cite[4.1.1]{BK} for  the representability \eqref{repcar}. The Cartier
dual can be described in the following way:
\begin{itemize}
\item For $M=[0\to G]$ we have $M^{\vee} = [T^{\vee}\by{u^{\vee}} A^{\vee}]$  where $T^{\vee}= \Hom (T,\G_m)$ is the character
group of $T$ and $u^{\vee}$ is the canonical homomorphism pushing out characters $T\to \G_m$ along the given extension $G$ of
$A$ by $T$.
\item For $M=[L\by{u} A]$ we have $M^{\vee} = [0\to G^u]$ where $G^u$ denote the group scheme which represents the functor
associated to $\Ext (M,\G_m)$. Here $\Ext (M,\G_m)$ consists of extensions of $A$ by $\G_m$ together with a trivialization of
the pull-back on $L$. In particular $[L\to 0]^{\vee} = \Hom (L,\G_m)$.
\item In general, the standard extension $M=[L\to G]$ of $M/W_{-2}M=[L\by{u} A]$ by $W_{-2}M = [ 0\to T]$ provides {\it via}
$\Ext (M/W_{-2}M,\G_m)$ the corresponding extension $G^{u}$ of $A^{\vee}$ by $\Hom (L,\G_m)$ and a boundary map
$$u^{\vee}:\Hom (W_{-2}M, \G_m)\to \Ext (M/W_{-2}M,\G_m)$$ lifting $T^{\vee}\to A^{\vee}$ as above.

\end{itemize}

\subsubsection{} A biextension is also providing natural pairings in realizations (see \cite[VIII 2]{SGA7} and 
\cite[10.2]{D}). In fact, for Grothendieck biextensions we also have an exact sequence 
$$\begin{array}{r} 0\to  \Ext (G_1\otimes G_2,H)\to \Biext (G_1,G_2;H) \to 
\Hom (\stor (G_1,G_2), H)\\[4pt]
\to  \Ext^2  (G_1\otimes G_2,H)\end{array}$$ and a natural map $T_{\ell}(G_1)\otimes T_{\ell}(G_2)\to T_{\ell}(\stor (G_1,G_2))$
(see \cite[VIII 2.1.13]{SGA7}) yielding a map
$$\Hom (\stor (G_1,G_2), H) \to \Hom  (T_{\ell}(G_1)\otimes T_{\ell}(G_2), T_{\ell}(H))$$ which in turns, by composition,
provides a map (see \cite[VIII 2.2.3]{SGA7})
$$\Biext (G_1,G_2;H) \to \Hom  (T_{\ell}(G_1)\otimes T_{\ell}(G_2), T_{\ell}(H)).$$

Similarly (non trivially! \cf \cite[\S 10.2.3-9]{D} and \cite{BD}) a biextension $P$ of 1-motives $M_1$ and $M_2$ by $H=\G_m$ provides the
following pairings: $T_{\ell}(M_1)\otimes T_{\ell}(M_2)\to T_{\ell}(\G_m)$ and $T_{DR}(M_1)\otimes T_{DR}(M_2)\to
T_{DR}(\G_m)$.  This latter pairing on De Rham realizations is obtained by pulling back $P$ to a {\it $\natural$-biextension}\,
$P^{\natural}$  of $M_1^{\natural}$ and
$M_2^{\natural}$ by $\G_m$. The Poincar\'e biextension $\cP_M$ of $M$ and $M^{\vee}$ by $\G_m$ is then providing compatibilities
between the Cartier dual of a 1-motive and the Cartier dual of its realizations. Moreover, over a base such that $p$ is locally
nilpotent, the Poincar\'e biextension is crystalline (see \cite[3.4]{AB}) providing  the {\it Poincar\'e crystal} of
biextensions $\cP_0^{\natural}$ of $M_0^{\natural}$ and 
$(M_0^{\vee})^{\natural}$ thus a pairing of $F$-crystals
$T_{crys}(M_0)\otimes T_{crys}(M_0^{\vee})\to T_{crys}(\G_m)$. We also have:
\subsubsection{\bf Scholium (\protect{\cite[10.2.3]{D}})}\label{Hodgepair}

{\em If $M_1$ and $M_2$ are defined over $\C$ then there is a natural isomorphism 
$$\Biext (M_1,M_2;\G_m)\cong \Hom_{\rm MHS}(T_{Hodge}(M_1)\otimes T_{Hodge}(M_2), \Z (1))$$} Over $\C$, all these pairings on
the realizations are deduced from Hodge theory.

\subsection{Symmetric avatar}\label{avatar} For a Deligne 1-motive $M=[L\by{u}G]$ and its Cartier dual $M^{\vee}=
[T^{\vee}\by{u^\vee}G^u]$ the Poincar\'e biextension $\cP_M= (\cP_A, \tau,\sigma)$ of $M$ and $M^{\vee}$ by $\G_m$ is
canonically trivialized on $L\times T^{\vee}$ by 
$\psi\df \tau_{\mid_{L\times T^{\vee}}} = \sigma_{\mid_{L\times T^{\vee}}}$  given by the push-out map $\psi_{\chi}: G\to\Chi
_*G$ along the character $\Chi : T \to \G_m$, \ie we have
$$\begin{array}{ccccccc} 0\to & T &\to & G &\to & A &\to 0\\ &\chi\downarrow\ &&\psi_\chi\downarrow\ \ &&\veq&\\ 0\to& \G_m&\to
& \Chi_*G&\to & A &\to 0\\
\end{array}$$ and $\psi (x,\chi) = \psi_{\chi}(u (x))\in \Chi _*G_{\bar u (x)}= (\cP_A)_{(\bar u (x), \bar u^{\vee} (\chi))}$
where $\bar u: L\to G \onto A$ and $ \Chi _*G = \bar u^{\vee} (\chi)$.  Actually, the data of $\bar u : L\to A$, $\bar u^{\vee}
: T^{\vee} \to A^{\vee}$ and $\psi$ determine both $M$ and $M^{\vee}$ under the slogan 
$${\rm trivializations}\ \ \iff \ \ {\rm liftings}$$ For example, for $\chi^1, \ldots ,\chi^r$ a basis of $T^{\vee}$ we can
regard $G$ as the the pull-back of $A$ diagonally embedded in $A^r$ as follows
$$\begin{array}{ccccccc} 0\to & T &\to & G &\to & A &\to 0\\ &\veq&&\downarrow&&\downarrow&\\ 0\to& \G_m^{r}&\to
&\chi_*^1G\times \cdots\times \chi_*^rG &\to & A^r &\to 0\\
\end{array}$$ and $(\psi (x,\chi^1), \ldots , \psi (x,\chi^r))$ provides a point of $G$ lifting $\bar u (x)$.
\subsubsection{} The {\em symmetric avatar}\, can be abstractly defined as $(L\by{u} A, L'\by{u'} A', \psi)$  where $L$, $L'$
are lattices, $A'$ is dual to $A$ and $\psi :L\times L'\to (u\times u')^*(\cP_A)$ is a trivialization of the Poincar\'e
biextension when restricted to $L\times L'$ (\cf \cite[10.2.12]{D}). In order to make up a category we define morphisms between
symmetric avatars by pairs of commutative squares such that the trivializations are compatible, \ie a map $$(L_1\by{u_1} A_1,
L_1'\by{u_1'} A_1', \psi_1)\to (L_2\by{u_2} A_2, L_2'\by{u_2'} A_2', \psi_2)$$ is a map $f : A_1\to A_2$ along with its dual $f
': A_2'\to A_1'$ and a pair of liftings  $g : L_1\to L_2$ of $fu_1$ and  $g' : L_2'\to L_1'$ of $f'u_2'$ such that 
$$\psi_1\mid_{L_1\times L_2'}= \psi_2\mid_{L_1\times L_2'}.$$ Here we have used the property that $(f\times 1)^*(\cP_{A_2}) =
(1\times f')^*(\cP_{A_1})$ for Poincar\'e biextensions. Denote $\M^{\rm sym}$ this category.

\subsubsection{\bf Scholium (\protect{\cite[10.2.14]{D}})}\label{equiavatar} {\em There is an equivalence of categories
$$M\mapsto (L\by{\bar u} A, T^{\vee}\by{\bar u^{\vee}} A^{\vee}, \psi): \M^{\fr}\longby{\simeq} \M^{\rm sym}$$ Under this
equivalence Cartier duality is 
$$(L\by{\bar u} A, T^{\vee}\by{\bar u^{\vee}} A^{\vee}, \psi)\mapsto (T^{\vee}\by{\bar u^{\vee}} A^{\vee}, L\by{\bar u} A,
\psi^t).$$} 
\subsubsection{} For the sake of exposition we sketch how to construct a map of symmetric avatars out of any biextension (almost
proving \eqref{repcar}, see \cite[4.1]{BK} for more details). Let $(P, \tau, \sigma)$ be a  $\G_m$-biextension of Deligne
1-motives $M_1$ and $M_2$. Translating {\em via}\, extensions, $P$ corresponds to a map $f :A_1\to A_2^{\vee}$ and the 
trivialization $\tau$ corresponds to a lifting $g: L_1 \to T_2^{\vee}$ of $f\bar u_1 $. Here we have $P\leadsto f\bar u_1
\leadsto 0\in \Hom (L_1,\Ext (G_2, \G_m))=  \Biext (L_1,G_2;\G_m)$ where
$$0\to \Hom (G_2,\G_m) \to T_2^{\vee}\by{\bar u_2^{\vee}} A_2^{\vee} \to \Ext (G_2, \G_m)\to 0$$ granting the existence of $g$
such that $\bar u_2^{\vee}g = f\bar u_1 $. Moreover $f\bar u_1 \leadsto [E]=0\in \Ext (L_1\otimes G_2,\G_m)=\Biext
(L_1,G_2;\G_m)$ and any section (= trivialization) of such trivial $\G_m$-extension $E$ is exactly given by an element $g\in 
\Hom (L_1, T_2^{\vee})=\Hom (L_1\otimes T_2, \G_m)$ as above. Since $P^t$ corresponds to the dual  $f^{\vee}: A_2\to A_1^{\vee}$
we have that 
$\sigma$ also corresponds to a lifting $g': L_2 \to T_1^{\vee}$ of $f^{\vee}\bar u_2 $  yielding $\bar u_1^{\vee}g' =
f^{\vee}\bar u_2$. Moreover, since  $P$  is a pull-back of $(f\times 1_{A_2})^{*}(\cP_{A_2}^t)$ then the trivialization $\tau$ is
the pull-back along $g\times 1: L_1\times G_2\to T_2^{\vee}\times G_2$ of the canonical trivialization $\psi_2^t$ on
$T_2^{\vee}\times G_2$ given by the identity.\footnote{Note that for the Poincar\'e biextension we have that the resulting $f$,
$f^{\vee}$, $g$ and $g'$ are all identities.} Since $P$ is also a pull-back of $(1_{A_1}\times f^{\vee})^{*}(\cP_{A_1})$ the
trivialization $\sigma$ on $G_1\times L_2$ is the pull-back of the canonical trivialization $\psi_1$ on $G_1\times T_1^{\vee}$
along 
$1\times g': G_1\times L_2\to G_1\times T_1^{\vee}$. Thus, if we further pull-back to $L_1\times L_2$ we get
$$\psi_2^t\mid_{L_1\times L_2}=\tau\mid_{L_1\times L_2}= \sigma\mid_{L_1\times L_2}=\psi_1\mid_{L_1\times L_2}$$ by assumption.
We therefore get a map
 $$(L_1\by{\bar u_1} A_1, T_1^{\vee}\by{\bar u_1^{\vee}} A_1^{\vee}, \psi_1)\to 
 ( T_2^{\vee}\by{\bar u_2^{\vee}} A_2^{\vee}, L_2\by{\bar u_2} A_2, \psi_2^t)$$

which turns back a map $M_1\to M_2^{\vee}$.

\subsection{1-motives with torsion} \label{torsion} An {\em effective}\, 1-motive which admits torsion (see \cite[\S 1]{BRS} and \cite{BK}) is $M= [ L \by{u} G]$ where $L$ is a locally constant (for the \'etale  topology) finitely generated abelian group
and $G$ is a semi-abelian scheme. Here
$L$ can be represented by an extension 
$$0\to L_{\tor}\to L \to L_{\fr}\to 0$$ where $L_{\tor}$ is finite and $L_{\fr}$ is free. An {\it effective} map from $M = [ L
\by{u} G]$ to $M'= [ L' \by{u'} G']$ is a commutative square and $\Hom_{\eff}(M, M')$ denote the abelian group of effective 
morphisms. The corresponding category is denoted
$\M^{\eff}$. We clearly have that $\M^{\fr}\subset \M^{\eff}$. For  $M = [ L \by{u} G]$ we set (see \cite[\S 1]{BRS} and
\cite{BK})
\begin{equation}\label{torfr}
\begin{array}{l} M_{\fr}\df [L_{\fr}  \by{\bar u}  G/u(L_{\tor})]\\[8pt] M_{\tor}\df [\ker (u)\cap L_{\tor}\to 0]\\[8pt]
M_{\tf}\df [L/\ker (u)\cap L_{\tor}  \by{u} G]\\[4pt]
\end{array}
\end{equation} considered as effective 1-motives. We say that $M$ is {\em torsion}\, if $L$ is torsion and $G=0$, $M$ is {\em
torsion-free}\, if  $\ker (u)\cap L_{\tor}=0$ and {\em free}\, if $L$ is free. There  are canonical effective maps $M \to
M_{\tf}$, $M_{\tor}\to M$ and $M_{\tf} \to M_{\fr}$.
\subsubsection{} A {\em quasi-isomorphism}\, (\qi for short) of 1-motives $M \to M'$ is a \qi of complexes of group schemes.
Actually, an effective map of 1-motives $M =[L\by{u}G]\to M'=[L'\by{u'}G']$ is a \qi of complexes if and only if we have that
$ \ker (u)=\ker (u')$ and $ \coker (u)=\coker (u')$ and thus $\ker$ and $\coker$ of $L\to L'$ and $G\to G'$ are equal. Then
$\coker (G\to G')= 0$, since it is connected and discrete, and $\ker (G\to G')$ is a finite group. Therefore a \qi of 1-motives
is given by  an isogeny $G\to G'$ such that $L$ is the pull-back of $L'$, \ie 
$$\begin{array}{ccccccc} 0\to & E &\to & G &\to & G' &\to 0\\ &\veq\ &&\ssp{u}\uparrow\ \ &&\ssp{u'}\uparrow\ &\\ 0\to & E &\to
& L &\to & L' &\to 0\\
\end{array}$$  where $E$ is a finite group.  We then define morphisms of 1-motives by localizing $\eM$ at the class of \qi and
thus set 
$$\Hom (M, M') \df \limdir{\qi} \Hom_{\eff} (\tilde M, M')$$ where the limit is taken over \qi $\tilde M \to M$ as above. We
then have a well-defined composition of morphisms of 1-motives (see \cite[1.2]{BRS})
$$\Hom(M,M') \times \Hom(M',M'') \to \Hom(M,M'').$$ In fact, for any effective morphism
$\tilde M \to M' $ and any \qi
$ \tilde M' \to M' $, there exists a \qi
$ \hat M \to \tilde M $ together with an effective morphism
$\hat M \to \tilde M' $ making up a commutative diagram (such that
$ \hat M \to \tilde M' $ is uniquely determined).
\subsubsection{} Denote the resulting category by $\M$, \ie objects are effective 1-motives and morphisms from $M$ to $M'$ can
be represented by a \qi $\tilde M \to M$ and an effective morphism $\tilde M \to M'$. This category has been introduced in
\cite{BRS} and it is further investigated in \cite{BK}. The main basic facts are the following:
\begin{itemize}
\item $\M$ is an abelian category where exact sequences can be represented by effective exact sequences of two terms complexes;
\item $\M^{\fr} \subset \M$ is a Quillen exact sub-category such that $M\mapsto M_{\fr}$ is left-adjoint to the embedding, \ie
we have $\Hom_{\eff}(M_{\fr},M') = \Hom(M,M')$ for $M\in \M$ and $M'\in \M^{\fr}$.
\end{itemize} Actually, we have $$\Hom_{\eff}(M, M') = \Hom (M,M')$$ for $M\in \M$ and $M'\in \M^{\fr}$. Clearly, this is
according with a corresponding {\em Scholium}~\ref{freein}  for the functor $\iota: \M \to D^b(k_{\rm fppf})$ which is still
faithful but, in general, not full for effective morphisms. A key point in order to show that $\M$ is abelian is the following.
\subsubsection{\bf Scholium (\protect{\cite[Prop. 1.3]{BRS}})}\label{strict} {\em Any effective morphism
$M\to M'$ can be factored as follows
$$\begin{array}{c} M  \longby{} M'\\
\displaystyle
\searrow\hspace*{0.5cm} \nearrow\\
\tilde M
\end{array}$$ where $M\to \tilde M$ is an effective morphism such that the kernel of the morphism of semi-abelian varieties is
connected, \ie a \emph{strict  morphism}, and $\tilde M \to M'$ is a \qi}\\[0.5cm] For example,  in the following canonical
factorisation induced by \eqref{torfr}
$$\begin{array}{c} M  \longby{} M_{\fr}\\
\searrow\hspace*{0.4cm} \nearrow \\ M_{\tf}
\end{array}$$ the effective map $M \to M_{\tf}$ is a strict epimorphism with kernel
$M_{\tor}$ and $M_{\tf}\to M_{\fr}$ is a \qi  We then always have a canonical exact sequence in $\M$
$$0\to M_{\tor}\to M \to M_{\fr}\to 0$$ We further have that the Hodge realization (see \ref{Hodge}) naturally extends to $\M
(\C)$ (see \cite[Prop. 1.5]{BRS}) and the functor
$$T_{Hodge}: \M (\C)\longby{\simeq} {\rm MHS}_1$$ is an equivalence between the category of 1-motives with torsion over $\C$ and
the category of $\Z$-mixed Hodge structures introduced in \ref{Hodge} above. Similarly, the other realizations extend to $\M$,
\eg (\cf \cite{BRS}, \cite{BK} and \ref{flat}) let $M/\ell^{\nu}$ be the torsion 1-motive (= finite group) given by the cokernel
of $\ell^{\nu}: M\to M$ the effective multiplication by $\ell^{\nu}$ which is fitting in an exact sequence (of finite groups)
$$0 \to {}_{\ell^{\nu}}M \to {}_{\ell^{\nu}}L\to {}_{\ell^{\nu}}G \to M/\ell^{\nu} \to L/\ell^{\nu} \to 0$$
and set $$T_{\ell}(M)\df \liminv{\nu}  M/\ell^{\nu}$$ Remark that Cartier duality does not extends to $\M$: such category $\M$
is just an algebraic version of ${\rm MHS}_1$.
\subsubsection{ } In \cite{BK} (\cf \ref{isogeny} below) we also consider larger categories of {\it non-connected} 1-motives,
\eg $[L\to G]$ where $G$ is a reduced group scheme locally of finite type over $k$ such that $G^0$ is semi-abelian\footnote{Note
that this condition can also be achieved by Murre's axiomatic \cite[Appendix A.1]{BSAP}.} and $\pi_0(G)$ is finitely generated.
If $M=[L\to G]$ is non-connected we get an effective 1-motive
$$M^0\df [L^0\to G^0]$$ where $L^0\subseteq L$ is the subgroup of those elements mapping to $G^0$ and 
$$\pi_0(M)\df  [L/L_0\into \pi_0 (G)]$$ is a discrete object.
\subsection{1-motives up to isogenies} \label{isogeny} For any additive category $\cC$ denote $\cC^{\sQ}$ the
$\Q$-linear category obtained from $\cC$ by tensoring morphisms by
$\Q$.

Let $\cC_1\df C^{[-1,0]}(Shv(k_{\et}))$ be the category of complexes of \'etale sheaves of length $1$ over $\Spec k$. Then 
$\cC_1$ and $\cC_1^{\sQ}$ are  abelian categories. We may view $\M^{\fr}$ and $\M^{\eff}$ as full subcategories of
$\cC_1$, hence $\M^{\fr , \sQ }$ and $\M^{\eff  , \sQ }$  as a full subcategory of $\cC_1^{ \sQ}$. The abelian category of {\em
1-motives up to isogenies} can be regarded {\it via} the following equivalences 
$$\M^{\fr , \sQ}\cong \M^{\eff  , \sQ}\cong \M^{\sQ}$$
since torsion 1-motives vanish and \qi of 1-motives are isomorphism in $\M^{\eff  , \sQ}$. Furthermore, let $\M^{\rm nc}$ be the
full subcategory of $\cC_1$ consisting of non-connected 1-motives, \ie complexes of the form $[L\to G]$ where $L$ is finitely
generated and $G$ is a commutative algebraic group whose connected component of the  identity $G^0$ is semi-abelian (see
\cite{BK}). We have that $\M^{\rm nc}\subset \cC_1$  is an abelian (thick) subcategory of $\cC_1$. For $M\in \M^{\rm nc}$ we
have that $M^0\into M$ and $M^0\onto M^0_{\fr}$ are isomorphisms in $\M^{\rm nc, \sQ}$. Thus $\M^{\fr,\sQ}\longby{\simeq}
\M^{\rm nc ,\sQ}$ is an equivalence of abelian categories.   
\subsubsection{\bf Scholium (\protect{\cite[1.1.3]{BK}})}\label{isonc}
{\em The category of Deligne 1-motives up to isogeny is equivalent to the abelian
$\Q$-linear category given by complexes of \'etale sheaves $[L\to G]$ where $L$ is (locally constant) finitely generated and $G$
is a commutative algebraic group whose connected component of the identity $G^0$ is semi-abelian.\\ }

Finally, this category
$\M^{\sQ}$ is of cohomological dimension $\leq 1$, \ie if
$\Ext^i (M,M') \neq 0$, for $M, M'\in \M^{\sQ}$, then $i =0$ or $1$ (\cite[Prop. 3.13]{OR}) and, clearly, the {\em
Scholium}~\ref{freein} holds for $\M^{\sQ}$ as well.

\subsection{Universal realization and triangulated 1-motives} 
I briefly mention some results from \cite{V}, \cite{OR} and \cite{BK}. Considering  the derived category of Deligne 1-motives up to isogeny we have a `universal realization'  in Voevodsky's triangulated category of motives. Notably, this realization has a left adjoint: the `motivic Albanese complex'.

\subsubsection{ } Recall that any abelian group scheme may be regarded as an \'etale sheaf with transfers (see \cite{V} for this notion and \cf \cite{VL}). Moreover, a 1-motive $M =[L\to G]$ is a complex of \'etale sheaves where $L$ and the extension $G$ of $A$ by $T$ are clearly homotopy invariants. Thus a 1-motive $M$ gives rise to an effective complex of homotopy invariant \'etale sheaves with transfers, hence to an
object of $\DM_{-,\et}^{\eff}(k)$ (see \cite[Sect. 3]{V} for motivic
complexes over a field $k$). 

Regarding 1-motives up to isogeny Nisnevich sheaves will be enough as
$\DM_{-}^{\eff}(k;\Q) \cong \DM_{-,\et}^{\eff}(k;\Q)$
is an equivalence of triangulated categories (see \cite[Prop. 3.3.2]{V} and
\cite[Th. 14.22]{VL}).

The triangulated category of effective geometrical motives $\DM_{\gm}^{\eff}(k;\Q)$ is the full triangulated sub-category of $\DM_{-}^{\eff}(k;\Q)$ generated by motives of smooth varieties: here the motive of $X$ denoted $M (X)\in \DM_{-}^{\eff}(k)$ is defined in \cite{V} by the Suslin complex $C_*$ of the representable presheaf with
transfers $L (X)$ associated to $X$ smooth over $k$. The motivic complexes provided by 1-motives up to isogeny actually belong to $\DM_{\gm}^{\eff}(k;\Q)$ (\cf \cite{OR} and \cite{BK}).

\subsubsection{\bf Scholium (\protect{\cite[Sect. 3.4, on page 218]{V} \cite{OR}})} \label{Tot} {\em There is a fully faithful functor $$\Tot : D^b(\M^{\sQ}) \by{\simeq} d_{\leq 1}\DM_{\gm}^{\eff}(k)\subseteq \DM_{\gm}^{\eff}(k;\Q)$$ whose essential image is  the thick triangulated subcategory $d_{\leq 1}\DM_{\gm}^{\eff}(k)\subseteq \DM_-^{\eff}(k)$ generated by motives of smooth varieties
of dimension $\leq 1$.}\\

Actually, in \cite{BK} we show that $D^b(\M^{\fr}) = D^b(\M)$ and we also refine this embedding to an integrally defined embedding of $ D^b(\M)[1/p]$ (where $p$ is the exponential characteristic) into the \'etale version $\DM_{\gm, \et}^{\eff}(k)$ of Voevodsky's category. The homotopy t-structure on $\DM_{-, \et}^{\eff}(k)$
induces a t-structure on $ D^b(\M)\cong d_{\leq 1}\DM_{\gm, \et}^{\eff}(k)$ with
heart the category $Shv_1 (k)$ of 1-motivic sheaves.
Here we also have that $\Tot ([0\to \G_m]) = \G_m[-1]\cong \Z (1)$ (see \cite[Th. 4.1]{VL} and \cite{BK}).

\subsubsection{ }  For $M\in\DM_{\gm}^{\eff}$ there is an internal  (effective) $\ihom (M, -)\in  \DM_{-}^{\eff}$ (see \cite[3.2.8]{V}). Set 
\begin{equation}\label{D1}
D\1 (M) \df\ihom (M, \Z (1))
\end{equation}
for any object $M\in \DM_{\gm}^{\eff}$. Actually, $D\1 (M)\in  d_{\leq 1}\DM_{\gm}^{\eff}$ (see \cite[3.1.1]{BK}) and restricted to $d_{\leq 1}\DM_{\gm}^{\eff}$ is an involution (see \cite[3.1.2]{BK}).

On the other hand, Cartier duality for 1-motives $M \mapsto M^{\vee}$ is an exact functor and extends to $D^b(\M)$. A key ingredient of \cite{BK} is that, under $\Tot$, Cartier duality is transformed into the involution
$M\mapsto \ihom (M,\Z (1))$ on $d\1\DM^{\eff}_{\gm}(k;\Q)$ given by
the internal (effective) Hom above.

\subsubsection{\bf Scholium (\protect{\cite[4.2]{BK}})} \label{Car} {\em We have a natural equivalence of functors
$$ \eta :  ( \ \ )^{\vee} \longby{\simeq}  \Tot^{-1} D\1\Tot $$
\ie under the equivalence $\Tot$ we have
$$\begin{array}{ccc}
 D^b(\M^{\sQ}) & \longby{\simeq} & d\1\DM_{\gm}^{\eff}(k;\Q)\\
\mbox{\rm \tiny ( \ )}^{\vee}\downarrow&  &\downarrow \mbox{\tiny D}\1\\
 D^b(\M^{\sQ}) & \longby{\simeq} & d\1\DM_{\gm}^{\eff}(k;\Q)
\end{array}$$}

Regarding $\Tot$ as the universal realization functor we expect that any other realization of $D^b(\M^{\sQ})$ (hence of $\M^{\sQ}$) will be obtained from a realization of $\DM_{\gm}^{\eff}(k;\Q)$ by composition with $\Tot$ and Cartier duality will be interchanging homological into cohomological theories.

\subsubsection{ } We show in \cite{BK} that $\Tot$ has a left adjoint
$\LAlb$. Dually,
composing with Cartier duality, we obtain $\RPic$. In order to construct $\LAlb$, let
$$d\1\df D_{\le 1}^2: \DM_{\gm}^{\eff}(k;\Q)
\to d\1\DM_{\gm}^{\eff}(k;\Q)$$
denote the functor
\begin{equation}\label{D2}
d\1(M)=\ihom (\ihom (M, \Z (1)), \Z (1))\in d\1\DM_{\gm}^{\eff}(k;\Q).
\end{equation}
The evaluation map yields a canonical map $a_M :M \to d\1 (M)$
that induces an isomorphism
$$\Hom (d\1 M, M')\iso \Hom (M, M')$$
for $M\in \DM_{\gm}^{\eff}(k;\Q)$ and $M'\in
d\1\DM_{\gm}^{\eff}(k;\Q)$. In fact, $M' =D\1 (N)$ for some $N\in d\1\DM_{\gm}^{\eff}(k;\Q)$ and if $C$ is the cone of $a_M$ then
\begin{multline*}
\Hom (C, M') = \Hom (C, D\1 (N)) = \Hom (C, \ihom (N,
\Z(1)))\\ =\Hom (C\otimes N, \Z (1)) = \Hom (N\otimes C, \Z
(1)) =\Hom (N, D\1 (C)) =0
\end{multline*}
since $D_{\le 1}^3=D_{\le 1}$.

\subsubsection{\bf Scholium (\protect{\cite[Sect. 2.2]{BK}})} \label{LAlb} {\em 
Define 
$$\LAlb : \DM_{\gm}^{\eff}(k;\Q)\to D^b(\M^{\sQ})$$
as the composition of $d\1\df D\1^2$ in \eqref{D2} and $\Tot^{-1}$. It is left adjoint to the embedding
$$\Tot : D^b(\M^{\sQ})\into \DM_{\gm}^{\eff}(k;\Q)$$ and $M\mapsto
a_M$ is the unit of this adjunction.}\\

The Cartier dual of $\LAlb$ is $\RPic = \Tot^{-1} D\1$.

\subsubsection{ }\label{1-motLR}
These functors provide natural complexes of 1-motives (up to isogeny)
 of any algebraic variety $X$ over a field $k$ if $\car (k) =0$ (for $X$ smooth and $k$ perfect even if $\car (k) >0$).  Their basic properties are investigated in \cite{BK}. We have:
\begin{itemize}
\item $\LAlb(X)\df \LAlb(M (X))$ the {\it homological Albanese complex}\, which is covariant on $X$ and, \eg it is homotopy invariant and satisfies Mayer-Vietoris;
\item $\LAlb^c (X)\df \LAlb(M^c (X))$ the {\it Borel-Moore Albanese complex}\, which is covariant for proper morphisms and  $$\LAlb(X) = \LAlb^c(X)$$ if $X$ is proper;
\item $\LAlb^*(X)\df \LAlb (M (X)^*(n)[2n])$ the {\it cohomological Albanese complex}\,
of $X$ purely $n$-dimensional, which is contravariant for maps between varieties of the same dimension and
$$ \LAlb^c (X) =\LAlb^*(X) $$
if $X$ is smooth (by motivic Poincar\'e duality $M^c(X)= M (X)^*(n)[2n]$, see \cite[Th.
4.3.2]{V});
\end{itemize}
and the Cartier duals:
\begin{itemize}
\item $\RPic(X)\df \RPic (M (X))$  the {\it cohomological Picard complex}\, which is contravariant in $X$;
\item $\RPic^c(X)\df \RPic (M^c (X))$ the {\it compactly supported Picard complex}\
such that 
$$\RPic(X) = \RPic^c(X)$$ if $X$ is proper;
\item $\RPic^*(X)\df \RPic (M (X)^*(n)[2n])$ the {\it homological Picard complex}\, of $X$ purely $n$-dimensional, which is covariant for maps between varieties of the same dimension and
$$\RPic^c(X) = \RPic^*(X)$$ 
if $X$ is smooth.\\
\end{itemize} 
Remark that the unit 
\begin{equation}\label{amap} 
a_X: M (X)\to \Tot \LAlb(X)
\end{equation} provide a universal map in
$\DM^{\eff}_{\gm}(k;\Q)$, the motivic Albanese map, which is an
isomorphism if $\dim (X)\leq 1$ and it refines the classical Albanese map and the less classical map in \cite{SZ}.

\subsection{1-motives with additive factors}  In order to keep care of non homotopical invariant theories we do have to include
additive factors. This is also suitable in order to include, in the 1-motivic world, the universal $\G_a$-extension
$M^{\natural}$ of a Deligne 1-motive $M$. In order to make Cartier duality working we cannot simply take $[L\to G]$ where $L$ is
(free) finitely generated and $G$ is a (connected) algebraic group: the Cartier dual of $\G_a$ is the formal group $\hat{\G}_a$,
\ie the connected formal additive $k$-group (see \cite{FO} and \cite{SGA3} for formal groups). Laumon \cite{LAU}  introduced a
generalization of Deligne's 1-motives in the following sense.   
\subsubsection{ } Laumon's 1-motives over a field $k$ of characteristic zero are given by $$M \df [F \by{u} G]$$ where $F$ is a
torsion free formal group and $G$  is a connected algebraic group, \ie $F$ has a presentation by a splitting extension
$$0\to F^0 \to F\to F_{\et}\to 0$$ where $F_{\et}$ \'etale over $k$ is further assumed torsion free (which means $F_{\et}(\bar
k)=\Z^r$) and $F^0$ is infinitesimal (that is given by a finite number of copies of $\hat{\G}_a$) and $G$ has a presentation
$$0\to T + V \to G\to A\to 0$$ where $T$ is a $k$-torus, $V$ is a $k$-vector group and $A$ is an abelian  variety. The map $u:
F\to G$ is any map of abelian fppf-sheaves so that an effective map $M\to M'$ is given by a map of complexes concentrated in
degrees $-1$ and $0$. Let $\M^{a,\fr}$ denote this category. 

\subsubsection{ } Recall \cite[2.2.2]{SGA3} that we have an antiequivalence  between (affine) algebraic groups and (commutative) formal groups, and, moreover, the  following formula (see \cite[5.2.1]{LAU}) holds: if such a formal group $F$ has Cartier dual 
$F^{\vee}$ and $A$ has dual $\Pic^0 (A) = A^{\vee}$ then
$$\Hom (F, A) = \Ext (A^{\vee},F^{\vee}).$$ Note that if $F = F^0$ is infinitesimal then $F^{\vee}\df \Lie (F)^{\vee}$ (= dual
$k$-vector space of the Lie algebra) and the extension associated to $F\to A$ is here obtained from the  universal
$\G_a$-extension $\Pic^{\natural}$ of $A^{\vee}$ by push-out  along $H^0(A, \Omega^1_A)=\Lie (A)^{\vee}\to \Lie (F)^{\vee}$.
The Cartier dual (\cf \ref{Cartier}) of $M =  [F \by{u} G]$ is given by an extension $G^u$ of $A^{\vee}$ by $F^{\vee}$
associated to the composite $F\to G\onto A$ and a lifting of $u^{\vee} : (T+V)^{\vee}\to A^{\vee}$ to  $G^u$ yielding
$$M^{\vee}\df [(T+V)^{\vee}\longby{u^{\vee}} G^u]$$ Moreover the Poincar\'e biextension of $M$ and $M^{\vee}$ by $\G_m$  is
obtained by pull-back from that of $A$ and $A^{\vee}$ as usual (see \cite[5.2]{LAU} for details).

\subsubsection{ } We have the following paradigmatic examples (\cf \ref{para} and 
\cite[5.2.5]{LAU}). If $X$ is a proper $k$-scheme then $[0\to \Pic^{0}_{X/k}]$ is a 1-motive defined by the Picard functor whose
Cartier dual (= the homological Albanese 1-motive) is  $[F\to \Alb (X)]$ where $\Alb (X) = \coker (\Alb (X_1)\to \Alb (X_0))$ is
dual to the abelian quotient of $\Pic^{0}_{X/k}$, $F_{\et}=\Z^r$ is the character group of the torus, see \eqref{semisimp}, and
$F^0=\hat{\G}_a^d$ corresponds to $d$-copies of $\G_a$ in $\Pic^{0}_{X/k}$.  Let $A$ be an abelian variety and let $[0\to
\Pic^{\natural ,0}_{A/k}]$ the 1-motive determined by the universal $\G_a$-extension of the dual $A^{\vee}$. The Cartier dual is
$[\hat{A}\to A]$ where $\hat{A}$ is the the completion at the origin of $A$. For $A= \Alb (X)$ and $X$ smooth proper over $k$
(of zero characteristic) we have so described the Cartier dual of $[0\to \Pic^{\natural ,0}_{X/k}]$.

\subsubsection{ } It seems possible to modify such a category, as we did (see \ref{torsion}) for Deligne 1-motives, in order to
include torsion, obtaining an abelian category. Just consider effective 1-motives $M = [F \to G]$ where $F$ is any formal group,
so that $F_{\et}$ may have torsion. However Cartier duality doesn't extend (here $F^{\vee}$ would be any, also non connected,
algebraic group). Let $\M^a$ denote this category. Similarly (\cf {\em Scholium}~\ref{isonc}) the category of Laumon 1-motives
up to isogeny is equivalent to the abelian $\Q$-linear category given by complexes of sheaves $[F\to G]$ where $F$ is a formal
group and $G$ is a commutative algebraic group. 

\subsubsection{ } A related matter is the Hodge theoretic counterpart of Laumon's 1-motives over $\C$ providing a generalized Hodge structure catching such additive factors (see \cite{FHS}).

 Provisionally define a {\it formal Hodge structure\/} (of level $\leq 1$) as follows. A formal group $H$ and a two steps filtration on a $\C$-vector space $V$, \ie $H=H^0\times H_{\et}$,  $H_{\et}=\Z^r + {\rm torsion}, H^0=\hat{\C}^d$ and $V^0\subseteq V^1\subseteq  V=\C^n$, along with a
mixed Hodge structure on the \'etale part, \ie say $H_{\et}\in {\rm MHS}_1$ for short, and a map $v: H\to V$. Regarding the induced map $v_{\et}: H_{\et} \to V$ we require the following conditions: for $H_{\sC}\df H_{\et} \otimes \C$ with Hodge filtration $F^0_{Hodge}$ and $c : H_{\et} \to H_{\sC}/F^0_{Hodge}$ the canonical map, the following
\begin{equation}\label{cond}
\begin{array}{ccc}
H_{\et}&\longby{v _{\et}}&V\\
\mbox{\tiny c}\downarrow & &\downarrow \mbox{\tiny pr}\\
H_{\sC}/F^0_{Hodge}&\longby{\simeq}& V/V^0
\end{array}
\end{equation}
commutes in such a way that $v _{\et}$ yields an isomorphism 
$H_{\sC}/F^0_{Hodge} \cong V/V^0$ restricting to an isomorphism
$W_{-2}H_{\sC}\cong V^1/V^0$.

Denote $(H, V)$ for short such a structure and let ${\rm FHS}_1$ denote the category whose objects are $(H, V)$ and the (obvious) morphisms given by commutative squares
compatibly with the data and preserving the conditions \eqref{cond}, \eg inducing a map of mixed Hodge structures on the \'etale parts.
Here we then get a forgetful functor $(H, V)\mapsto H_{\et}$ from  ${\rm FHS}_1$ to  ${\rm MHS}_1$, left inverse of the embedding $H\mapsto (H, H_{\sC}/F^0_{Hodge})$. Actually we can define $(H, V)_{\et}\df (H_{\et}, V/V^0)$ and say that a formal Hodge structure is {\it \'etale}\,  if $(H, V) = (H, V)_{\et}$, \ie if $H^0=V^0=0$.  The full subcategory ${\rm FHS}_1^{\et}$ of \'etale structures is then equivalent to ${\rm MHS}_1$ {\it via} the forgetful functor and the functor $(H, V)\mapsto (H, V)_{\et}$ is a left inverse of the 
inclusion ${\rm FHS}_1^{\et}\subset {\rm FHS}_1$. Remark that $(H,V)$ with $H_{\et}$ pure of weight zero exists if and only if $V=V^1=V^0$. Thus if $v$ restricts to a map $v^0: H^0\to V^0$ then $(H^0,V^0)$ is a formal substructure of $(H, V)$ and we have
a `non canonical' extension
\begin{equation}\label{formhodge} 0\to (H^0, V^0)\to (H,V) \to (H,V)_{\et}\to 0
\end{equation} 
Say that $(H, V)$ is {\it connected}\, if $(H, V)_{\et}= 0$ and that it is {\it special}\, if $(H^0, V^0)\df  (H, V)^0$ is a substructure of $(H,V)$ or, equivalently,  $(H, V)_{\et}$ is a quotient of $(H,V)$: the above extension \eqref{formhodge} is then characterizing  special structures. 

\subsubsection{ } Extending Deligne's Hodge realization (\cf \ref{Hodge}) for a given 1-motive $M = [F \to G]$ consider the pull-back $T_{\oint}(F)$ of $F \to G$ along $\Lie (G)\to G$. Here $T_{\oint}(F)$  is a formal group and the canonical map $T_{\oint}(F)\to \Lie (G)$ provides the `formal Hodge realization' of $M$
$$T _{\oint}(M)\df (T _{\oint}(F), \Lie (G))$$
as follows.

For $M = [F \by{u} G]$ over $k$ let $V (G)\df \G_a^n\subseteq G$ be the additive factor and display $G$ as follows
\begin{equation}\label{gmext} 0\to V (G)\to G \to G_{\times}\to 0
\end{equation} where $G_{\times}$ is the semi-abelian quotient. We have that $\Lie (G)$ is the pull-back of $\Lie (G_{\times})$ and
$H_1(G)=H_1(G_{\times})$. Moreover $F=F^0\times_{k}F_{\et}$ (canonically) and we can set $M_{\et}\df [F_{\et}\to G_{\times}]$. The functor $M\mapsto M_{\et}$ is a left inverse of the inclusion of Deligne's 1-motives. We have
that $T_{\oint}(F)_{\et}$ is an extension of $F_{\et}$ by $H_1(G_{\times})$ so that, by construction, the formal group $T_{\oint}(F)$
has canonical extension 
$$0\to F^0\to T _{\oint}(F) \to T_{\sZ}(M_{\et})\to 0$$ where $ T_{\sZ}(M_{\et})$ is the $\Z$-module of the usual Hodge
realization (see \ref{Hodge}) providing the formula $T _{\oint}(F)_{\et} = T_{\sZ}(M_{\et})=$ the pullback of $F_{\et}\into F$
along $T _{\oint}(F)\to F$.

Thus $(T _{\oint}(F), \Lie (G))\in {\rm FHS}_1$ where $T_{\oint}(F)_{\et}$ is the underlying group of the Hodge structure $T_{Hodge}(M_{\et})$, the filtration $V(G)\subseteq V(G) + \Lie (T)\subseteq \Lie (G)$ is the two steps filtration and the condition \eqref{cond} is provided by construction (see \ref{Hodge}) since $T_{\sC}(M_{\et})\df T_{\sZ}(M_{\et})\otimes \C \cong \Lie (G^{\natural})$ (see \ref{DeRham}, here $M_{\et}^{\natural} =[F_{\et}\to G^{\natural}]$  is the universal $\G_a$-extension of $M_{\et}$), \ie 
$$
\begin{array}{ccc}
T_{\sZ}(M_{\et})&\longby{v _{\et}}&\Lie (G)\\
\mbox{\tiny c}\downarrow & &\downarrow \mbox{\tiny pr}\\
T_{\sC}(M_{\et})/F^0_{Hodge}&\longby{\simeq}&\Lie (G_{\times})
\end{array}
$$
commutes and $W_{-2}T_{\sC}(M_{\et})\cong \Lie (T)$.

Moreover, if $u$ restricted to $F^0$ is mapped to $V(G)$ we can further set $V (M)\df [F^0\to V(G)]$ fitting in an extension
\begin{equation}\label{formext} 0\to V (M)\to M \to M_{\et}\to 0
\end{equation} providing a `non canonical' extension of the 1-motive (\cf \eqref{formhodge}, here \eqref{formext} becomes \eqref{formhodge} by applying $T _{\oint}$). Note that $M_{\et}$ is pure of weight zero if and only if $T _{\oint}(M) = M$. 
For example, if $W\by{u} V$ is a linear map between $\C$-vector spaces, and $M = [\hat{W}\by{u} V]$ is the  induced 1-motive (here $\hat{W}$ is the formal completion at the origin, \cf \cite[5.2.5]{LAU}) then $T _{\oint}(M) = M$.

\subsubsection{\bf Scholium (\protect{\cite{FHS}})}\label{form} {\em There is an equivalence of categories
$$M\leadsto T _{\oint}(M): \M^{a,\fr}(\C)\longby{\simeq} {\rm FHS}_1^{\fr}$$
between Laumon's 1-motives and torsion free formal Hodge structures (of level$\leq 1$)
providing a diagram
$$\begin{array}{ccc}
\M^{\fr}(\C)&\by{\simeq}& \mbox{\rm MHS}_1^{\fr}\\
\uparrow\downarrow & &\uparrow\downarrow\\
\M^{a,\fr}(\C)&\by{\simeq}& \mbox{\rm FHS}_1^{\fr}
\end{array}$$}

Regarding duality for $(H, V)= T _{\oint}(M)$ such that $H_{\et}$ is free, we can
argue cheaply defining it as follows
$$T _{\oint}(M)^{\vee}\df T _{\oint}(M^{\vee})$$
After Cartier duality (\cf \cite[5.2]{LAU}) it is easy to check that this is a self duality extending the one on ${\rm MHS}_1^{\fr}$. 
For example $H_{\et}^{\vee}\df \ihom (H_{\et}, \Z (1)) = (H, V)^{\vee}_{\et}$ as usual and the Cartier dual of $(H,V)$ with $H_{\et}$ of weight zero such that \eqref{formext} splits  $M= V(M)\oplus M_{\et}$, is obtained as follows: $V(M)^{\vee}$ is given by
$\hat{V}^{\vee}\to \Lie (H^0)^{\vee}$ obtained from the induced map $\Lie (H^0)\to V$, taking the dual vector space map
$V^{\vee}\to \Lie (H^0)^{\vee}$ and its completion $\hat{V}^{\vee}\to V^{\vee}$ at the origin\footnote{Note that in 
characteristic zero there is a canonical equivalence of categories between  Lie algebras and infinitesimal formal groups.} thus
$$(H^0\times \Z (0)^{\oplus r}, V)^{\vee} = (\hat{V}^{\vee}\times \Z (1)^{\oplus r},\Lie (H^0)^{\vee}\times \G_a^{\oplus r})$$ where the map
$\Z (1)^{\oplus r}\to \G_a^{\oplus r}$ is canonically induced from the exponential map, \eg in particular $(\Z (0), 0)^{\vee} = (\Z (1), \G_a)$. 

\subsubsection{ }  Formal Hodge structures\footnote{We here mean to deal with arbitrary Hodge numbers. However, for the sake of brevity, no more details on FHS are provided: generalizing our definition above it's not that difficult but it's more appropriate to treat such a matter separately.} would pitch in the following diagram
$$\begin{array}{ccc}
\mbox{Deligne's 1-motives}&\longby{T _{Hodge}}& \mbox{MHS}\\
\uparrow\downarrow & &\uparrow\downarrow\\
\mbox{Laumon's 1-motives}&\longby{T _{\oint}}& \mbox{FHS}
\end{array}$$
where 
\begin{itemize} 
\item ${\rm FHS}$ would be a rigid tensor abelian category which is an enlargement of MHS and $H 
\mapsto H_{\et}$ would yield a functor from FHS to MHS, left inverse of the  embedding;
\item $T _{\oint}$ would be fully faithful so that under the
realizations Cartier  duality corresponds to a canonical
$\ihom (-,\Z_a (1))$ involution.\footnote{Here we clearly have the candidate $T _{\oint}([0\to \G_m]) \df \Z_a (1)$ and a formal version of {\em Scholium}\, \ref{Hodgepair} should be conceivable.}
\end{itemize}

\subsubsection{ } Similarly define other realizations, \eg see \cite{BAB} where  we obtain the {\em sharp}\, De Rham realization $T_{\sharp}$. For example, if $F^0 =0$ we can describe $T_{\sharp}$ out of the universal
$\G_a$-extension $M_{\et}^{\natural}$ (see \ref{DeRham}), defining the algebraic group $ G^{\sharp}$ by pull-back {\em via} \eqref{gmext}
as an extension 
$$0\to \Ext (M_{\et},\G_a)^{\vee}\to G^{\sharp}\to G\to 0$$  taking 
$\Lie (G^{\sharp})$. In this case we thus obtain a canonical extension
$$0\to V (G)\to M^{\sharp} \to M_{\et}^{\natural}\to 0$$ 
and we can relate to $(H, V)= T _{\oint}(M)$, where $H^0=F^0=0$, passing to Lie algebras, by the following pull-back diagram
$$\begin{array}{cccc}
&&0&\\
&&\uparrow&\\
0\to F^0_{Hodge}&\to H_{\sC}\to & H_{\sC}/F^0_{Hodge}&\to 0\\
\veq & \uparrow &\uparrow\\
0\to \Ext (M_{\et},\G_a)^{\vee}&\to\Lie (G^{\sharp}) \to& \Lie (G)&\to 0\\
& &\uparrow&\\
& &V(G)& \\
& &\uparrow&\\
& &0&
\end{array}$$
where $H_{\et} = T_{\sZ}(M_{\et})$, $H_{\sC}=  \Lie (G^{\natural})$, $H_{\sC}/F^0_{Hodge} =  \Lie (G_{\times})$, $V^0=V(G) \subseteq V=\Lie (G)$.  Set $H_{\sC}^{\sharp}\df \Lie (G^{\sharp})$ and by the universal property we get
an induced map $c^{\sharp} : H_{\sZ} \to H_{\sC}^{\sharp}$ providing a splitting 
of the projection $H_{\sC}^{\sharp}\onto H_{\sC}$ and the diagram above can be translated in the following diagram
\begin{equation}\label{sharp}
\begin{array}{cccc}
H_{\sC} & \to& H_{\sC}/F^0_{Hodge}&\\
\uparrow &&\uparrow\\
H_{\sC}^{\sharp}&\to & V&\\
 \uparrow &&&\\
 H_{\sC} &&& 
\end{array}
 \end{equation}
\subsubsection{ } Remark that, from a different point of view Bloch and Srinivas \cite{BSE} proposed a category of enriched Hodge structures EHS whose
objects are pairs $E\df (H, V_{\d})$ where $H$ is a mixed  Hodge structure and $V_{\d}$ is a diagram (not a complex)
$\cdots = V_{a+1}=V_{a}\to V_{a-1}\to \cdots \to V_{b}\to 0\to 0 \cdots$ of $\C$-vector spaces such that $V_{i}\to
H_{\sC}/F^{i}$ (compatibly with  the diagram) and there is a map $H_{\sC}\to V_{a}$ such that $H_{\sC}\to  V_{a}\to
H_{\sC}/F^{a}$ is the identity, thus 
$F^a =0$. There is a canonical functor $E \mapsto H$ to MHS (with a right  adjoint). 
It is not difficult to see that ${\rm EHS}_1$ is equivalent to  ${\rm FHS}_1^s\subset  {\rm FHS}_1$ the subcategory of special formal Hodge structures given by \eqref{formhodge}.
Sharp De Rham realization also clearly provides an enriched Hodge structure, \eg {\it via}\, \eqref{sharp}. In fact, we can refine the construction  \eqref{sharp} obtaining a functor $$T_{\sharp}^s : {\rm FHS}_1^s\to {\rm EHS}_1$$ by sending
$$(H, V)\mapsto T_{\sharp}^s (H, V) \df (H_{\et}, H_{\sC}^{\sharp}\to V)$$ 
where $H_{\sC}^{\sharp}$ (along with the splitting) is just obtained by pull-back when $H^0 =0$ (see \cite{BAB} for the precise statements and further properties).

\section{On 1-motivic (co)homology}
In the previous section 2 we have provided realizations as covariant and contravariant functors from categories of 1-motives to categories of various kind of structures. Here we draft a picture (which goes back to the algebraic geometry constructions of section 1)  providing 1-motives, \ie 1-motivic cohomology, whose realizations are the `1-motivic part' of various existing (or forthcoming) homology and cohomology theories.

\subsection{Albanese and Picard 1-motives}
Let $X$ be a complex algebraic variety and let $H^*(X, \Z )$ be the mixed Hodge structure on the singular cohomology of the associated analytic space.
Denote $H^*_{(1)}(X, \Z (\cdot))\subseteq H^*(X, \Z (\cdot))$ the largest substructure and $H^*(X, \Z (\cdot))^{(1)}$ the largest quotient in ${\rm MHS}_1$ (\cf \ref{Hodge}).

\subsubsection{\bf Deligne's Conjecture (\protect{\cite[10.4]{D}})}\label{dconj} {\em
Let $X$ be a complex algebraic variety of dimension $\leq n$. There exist algebraically defined 1-motives whose Hodge realizations over $\C$ are
$H^i_{(1)}(X, \Z (1))_{\fr}$, $H^i(X, \Z (i))^{(1)}_{\fr}$ for $i\leq n$ and
$H^i(X, \Z (n))^{(1)}_{\fr}$ for $i\geq n$ and similarly for $\ell$-adic and De Rham realizations.}\\

The results contained in \cite{C}, \cite{BSAP}, \cite{RA}, \cite{RA1}, \cite{BRS} and \cite{AB} show some cases of this conjecture. Over a field $k$, $\car (k)=0$, with the notation of \cite{BSAP}:
\begin{itemize}
\item $\Pic^+(X)$ which reduces to \eqref{open} if $X$ is smooth  (or to the simplicial $\bPic^0$ if $X$ is proper) provides an algebraic definition of $H^1_{(1)}(X, \Z (1))_{\fr}=H^1(X, \Z (1))$;
\item $\Alb^+ (X)$ is an algebraic definition of $H^{2n-1}(X, \Z (n))_{\fr}$ for $n =\dim (X)$;
\item $\Pic^- (X)=\Alb^+ (X)^{\vee}$ is an algebraic definition of $H_{2n-1}(X, \Z (1- n))$;
\item $\Alb^- (X)=\Pic^+(X)^{\vee}$ which reduces to Serre Albanese if $X$ is smooth, is an algebraic definition of $H_{1}(X, \Z )_{\fr}$.
\end{itemize}  
Moreover, in \cite{BRS} we have constructed effective 1-motives with torsion:
\begin{itemize}
\item $\Pic^+(X,i)$ for $i\geq 0$ providing an algebraic definition of $H^{i+1}_{(1)}(X, \Z (1))_{\fr}$ up to isogeny.
\end{itemize} 
 
Actually, for $Y$ a closed subvariety of $X$ we have $\Pic^+(X,Y;i)$  (= $M_{i+1} (X, Y)$ in the notation of \cite{BRS}) such that $\Pic^+(X,\emptyset ;0)= \Pic^+(X)$. These $\Pic^+(X,Y;i)$ are obtained using appropriate `bounded resolutions' which also provide a canonical {\em integral}\, weight filtration $ W $ on the relative cohomology $H^{*}(X,Y;\Z)$ (see \cite[2.3]{BRS}).

\subsubsection{\bf Scholium (\protect{\cite[0.1]{BRS}})}\label{thm} {\em
There exists a canonical isomorphism of mixed Hodge structures
$$\phi_{\fr} : T_{Hodge}(\Pic^+(X,Y;i))_{\fr} \longby{\simeq}
W_{0}{H}_{(1)}^{i+1}(X,Y;\Z (1))_{\fr}$$
and similarly for the $l$-adic and de Rham realizations.}\\

This implies Deligne's conjecture on  ${H}_{(1)}^{*}(X,\Z (1))_{\fr}$ up to
isogeny and in cohomological degrees $\leq 2$ even without isogenies by dealing with such 1-motives with torsion. The conjecture without isogeny is reduced to
$${H}_{(1)}^{*}(X,Y;\Z)_{\fr} = W_{2}{H}_{(1)}^{*}(X,Y;\Z)_{\fr}.$$

Here the semiabelian part of $\Pic^+(X,Y;i)$ yields
$W_{-1}{H}_{(1)}^{j}(X,Y;\Z (1))_{\fr} $ and the torus corresponds to
$ W_{-2}{H}_{(1)}^{j}(X,Y;\Z (1))_{\fr} $.

\subsubsection{ } In \cite{AB} we have also formulated a corresponding statement \ref{dconj} for the crystalline realization. Recall that de Jong \cite[p. 51-52]{DJ} proposed a definition of crystalline cohomology\footnote{Note that we can
also deal with rigid cohomology and everything here can be rephrased
switching crystalline to rigid.} forcing cohomological descent. Let $X$ be an algebraic variety, over a perfect field $k$, de Jong's theory \cite{DJ} provide a pair $(X_{\d},
Y_{\d})$ where $X_{\d}$ is a smooth proper simplicial scheme,
$Y_{\d}$ is a normal crossing divisor in $X_{\d}$ and $X_{\d}-Y_{\d}$ is a smooth proper hypercovering of $X$.  Set $H^*_{\crys} (X/\W (k)) \df \HH^*_{\logcrys} (X_{\d} , {\rm Log\ } Y_{\d})$ where $(X_{\d}, {\rm Log\ } Y_{\d})$ here denotes the simplicial logarithmic structure on $X_{\d}$ determined by $Y_{\d}$ (see \cite[\S 6]{AB}). The question here 
(\cf \cite{DJ}) is that $H^*_{\crys} (X/\W (k))$ is not {\it a priori}\, well-defined. Similarly to \cite[2.3]{BRS} we may also expect a weight filtration $W_*$ on the crystalline cohomology $H^{*}_{\crys} ((X, Y)/\W (k))$ of a pair $(X, Y)$. 

Over a perfect field, using de Jong's resolutions, it is easy to obtain an appropriate construction of $\Pic^{+}(X,Y; i)$ such that $\Pic^{+}(X,\emptyset ; 0) = \Pic^{+}(X)$ as above. In \cite[Appendix A]{AB} we have shown that $\Pic^{+}(X)$ is really well-defined and independent of the choices of resolutions or compactifications. However, it is not clear,  for $i>0$, if $\Pic^{+}(X,Y; i)$  is {\it integrally} well-defined: the 1-motive is well-defined up to $p$-power isogenies in characteristic $p$ by \cite[A.1.1]{AB} and a variant of \cite[Thm. 3.4]{BRS} for $\ell$-adic realizations with $\ell\neq p$.

\subsubsection{\bf Crystalline Conjecture (\protect{\cite[Conj. C]{AB}})}\label{crysconj} {\em  Let $H^{*}_{\crys ,(1)} ((X, Y)/\W (k))$ denote the submodule of $W_2H^{*}_{\crys} ((X, Y)/\W (k))$
whose image in $\gr_W^2$ is generated by the image of the discrete
part of $\Pic ^{+}(X,Y; i)$ under a suitable cycle map. Then there
is a canonical isomorphism (eventually up to $p$-power isogenies)
$$T_{\crys} (\Pic ^{+}(X,Y; i)) \longby{\simeq} H^{i+1}_{\crys
,(1)} ((X, Y)/\W (k))(1)$$ of filtered $F$-$\W (k)$-modules (\ie we expect a crystalline analogue of \ref{dconj}).}\\

We can show this statement for $i =0$ and $Y =\emptyset $ (see \cite[Thm.~B${}^{\prime}$]{AB}).  The corresponding general statement for De Rham cohomology over a field of characteristic zero is \cite[Thm.~3.5]{BRS}. 

\subsubsection{ } According with the program in \cite{BK} these $\Pic ^{+}(X,Y; i)$ would get linked to Voevodsky's theory of triangulated motives as follows.
The covariant functor $M : Sm/k \to \DM_{\gm}^{\eff}(k)$ from the category of
smooth schemes of finite type over $k$ (a field admitting resolution of singularities) extends to all schemes of finite type (see \cite[\S 4.1]{V}). Thus the motivic Albanese complex $\LAlb (X)$ and the motivic Picard complex $\RPic (X)$ are well-defined for any such scheme $X$ (and similarly for the other complexes, see \ref{1-motLR}). Consider the (co)homology 1-motives (up to isogenies) $H_{i}(\LAlb (X))\df \LA{i}(X)$ and $H^i(\RPic (X))\df \RA{i}(X)$ for $i\in\Z$. We have that $$\RA{i}(X) = \LA{i}(X)^{\vee}$$ by motivic Cartier duality (see {\it Scholium}\, \ref{Car}).

\subsubsection{\bf LAlb - RPic Hypothesis (\protect{\cf \cite{BK}}).}\label{LRthm} {\em We assume the following picture (up to isogeny):\\
\begin{itemize}
\item $T_{Hodge}(\LA{i}(X)) = H_{i}(X,\Z)^{(1)}_{\fr}$ = 1-motivic singular homology mixed Hodge structure; 
\item $T_{Hodge}(\LA{i}^c(X)) = H_{i}^{BM}(X,\Z)^{(1)}_{\fr}$ = 1-motivic Borel-Moore homology mixed Hodge structure;
\item $T_{Hodge}(\LA{i}^*(X)) = H^{2n-i}(X,\Z (n))^{(1)}_{\fr}$ = 1-motivic Tate twisted singular cohomology mixed Hodge structure of $X$ $n$-dimensional;
\end{itemize}
and dually:
\begin{itemize}
\item $T_{Hodge}(\RA{i}(X)) = H^{i}_{(1)}(X,\Z (1))_{\fr}$ = 1-motivic singular cohomology mixed Hodge structure; 
\item $T_{Hodge}(\RA{i}^c(X)) = H^{i}_{c, (1)}(X,\Z (1))_{\fr}$  = 1-motivic compactly supported cohomology mixed Hodge structure; 
\item  $T_{Hodge}(\RA{i}^*(X)) = H_{2n-i, {(1)}}(X,\Z (1-n))_{\fr}$  = 1-motivic Tate twisted singular homology of $X$ $n$-dimensional.\\
\end{itemize} 
Similar statements for $\ell$-adic, De Rham and crystalline realizations are also workable  (providing a positive answer to \ref{dconj} and \ref{crysconj}).}\\
 
It is not difficult  (see \cite{BK}) to compute these 1-motivic (co)homologies for $X$ smooth or a singular curve. We recover in this way Deligne-Lichtenbaum motivic (co)homology of curves (\cf \cite{D} and \cite{LI}). The picture above also recover the previously mentioned Picard and Albanese 1-motives as follows
$$\LA{1}(X) = \Alb^-(X) \hspace*{2cm} \LA{1}^*(X) = \Alb^+ (X)$$
and 
$$\RA{1}(X) = \Pic^+(X) \hspace*{2cm} \RA{1}^*(X) = \Pic^- (X).$$
Finally, for $i\ge1$, we should get a formula like that
$$\RA{i}(X,Y) = \Pic ^{+}(X,Y; i-1)  \hspace*{2cm} \LA{i}(X,Y) =  \Alb ^{-}(X,Y; i-1)$$
with the obvious meaningful notation adopted above.

\subsection{Hodge 1-motives} We shortly explain the point of view developed in \cite{BA} extending Deligne's philosophy \ref{dconj} to algebraic cycles in higher codimension (\cf \ref{CH}). The starting point is by looking at the side of \ref{thm} which provides a Lefschetz theorem on $(1,1)$-classes, \ie in degrees $>1$. See also \cite{BRV}.

\subsubsection{ } \label{NS+}
 Define $\NS^{+} (X, Y; i)$ for $i\geq 0$ as the quotient of $\Pic ^{+}(X,Y; i+1)$ by its toric part and consider the extension 
$$0\to W_{-2} \to  \Pic ^{+}(X,Y; i+1)  \to \NS^{+} (X, Y; i)\to 0$$
It follows from \ref{thm} up to isogeny 
$$T_{Hodge}(\NS^{+} (X, Y; i))=W_0{H}_{(1)}^{2 + i}(X,Y;\Z (1))/W_{-2}$$
given by the extension 
$$0\to \gr^W_{1}\to W_2{H}^{2 + i}(X,Y;\Z )/W_{0}\to \gr^W_{2}\to 0$$
pulling back $(1,1)$-classes in $ \gr^W_{2}$ (and twisting by $\Z (1)$).
We may call  $\NS^{+} (X, Y; i)$ the {\it Hodge-Lefschetz} 1-motive since, \eg if $X$ is smooth proper and $Y=\emptyset$ we obtain $\NS^{+} (X; 0) = \NS (X)$ and $\NS^{+} (X; i) =0$ for $i\neq 0$. 

\subsubsection{ } Set $H\df {H}^{2p + i}(X,Y;\Z )$ for a fixed $p\geq 1$ and $i\geq -1$ and consider
 $$0\to\gr^W_{2p-1}H \to W_{2p}H/W_{2p-2}H\to \gr^W_{2p}H\to 0$$ 
 given by the integral weight filtration (see \cite{BRS}).
 Consider the integral $(p,p)$-classes $H^{p,p}_{\sZ}\df \Hom_{\rm MHS} (\Z (-p), \gr^W_{2p}H)$ and the associated intermediate jacobian  $J^p(H)\df \Ext (\Z (-p), \gr^W_{2p-1}H)$ which is just a complex torus if $p> 1$ (see \cite{CA}). Consider the largest abelian subvariety $A^p(H)$ of the torus $J^p(H)$ which corresponds to the maximal polarizable substructure of $gr^W_{2p-1}H$ purely of types $\{(p-1,p), (p, p-1)\}$.
Define the group of {\it Hodge cycles} $H^p(H)$ as the
preimage in $H^{p,p}_{\sZ}$ of $A^p(H)$ under the extension class map $e^p: H^{p,p}_{\sZ} \to J^p(H)$.  Define the {\it Hodge 1-motive} by
$$e^p : H^p(H) \to A^p (H)$$
and the corresponding mixed Hodge structure $H^h\in {\rm MHS}_1$.
\subsubsection{\bf Anodyne Hodge Conjecture (\cf \protect{\cite[2.3.4]{BA}})}\label{AHCconj} {\em Let $X$ be an algebraic variety and $Y$ a closed subvariety defined over a perfect field $k$.  There exist algebraically defined 1-motives with torsion $\Xi^{i,p}(X, Y)\in \M (k)$ whose Hodge realization over $k=\C$ are ${H}^{2p + i}(X,Y;\Z )^h\in {\rm MHS}_1$, \ie here ${H}^{2p + i}(X,Y;\Z )$ is the associated mixed Hodge structure (for $p\geq 1$ and $i\geq -1$) so that
$$T_{Hodge}(\Xi^{i,p}(X, Y)) \cong H^{2p+i}(X, Y; \Z)^h$$
 and similarly for $\ell$-adic, De Rham and crystalline realizations.}\\
 
For $p=1$ this is `almost' true ($\approx$ Deligne's conjecture \ref{dconj} and \ref{thm} but \ref{crysconj}) and it  follows from \ref{NS+} and $\Xi^{i,1}(X, Y)=\NS^{+} (X, Y; i)$. One can also easily formulate a homological version of \ref{AHCconj}. 
Recall that for $\bar X $ smooth proper purely $n$-dimensional and $Y + Z$ normal crossing divisors on $\bar X$  (in particular when $X = \bar X -Z$ and $Y\cap Z = \emptyset$) we have 
 $$ H^{2p+i}(\bar X -Z, Y; \Z (p)) \cong H_{2r -i}(\bar X - Y, Z;\Z (-r))\hspace*{0.7cm} (p=n-r)$$ as mixed Hodge structures (see \cite[2.4.2]{BSAP}). 
For $X$ smooth and proper (here we assume that $Y = Z = \emptyset$ and $X = \bar X$) we get $H^{2p+i}(X, \Z)^h\neq 0$ if and only if $i = -1, 0$ and \ref{AHCconj} reduces to the quest of an algebraic definition of $A^p\subseteq J^p$ or $H^{p, p}_{\sZ}$ respectively. Classical Grothendieck-Hodge conjecture then provides candidates up to isogeny.

\subsubsection{ } \label{univab}
For $X$ a smooth proper $\C$-scheme we can consider $J^p_a(X)\subseteq J^p(X)$ the image of $CH^p(X)_{\rm alg}$ (\cf \ref{CH}) under the Abel-Jacobi map: the usual Grothendieck-Hodge conjecture claims that $J^p_a(X)$ is the largest abelian variety in $ J^p(X)$, \ie that $A^p=J^p_a$ (up to isogeny) and ${H}^{2p -1}(X, \Z )^h$ is algebraically defined {\it via}\, the coniveau filtration. Similarly, the image of $\NS^p (X)$ generates $H^{p, p}_{\sZ}$ (with $\Q$-coefficients). In the most wonderful world (mathematics!?) the 1-motivic sheaf $(CH_{X}^p)^{(1)}$ in \ref{CH} could make the job providing an algebraically defined extension of $H^{p, p}_{\sZ}$ by $J^p_a$ compatibly with \eqref{NSext} (here $J^p_a$ would also coincide with the universal regular quotient of $CH^p(X)_{\rm alg}$ when $X$ is smooth and proper).

If $X$ is only proper then let $\pi : \Xs \to X$ be a resolution and consider the Chow groups of each component $X_i$ of $\Xs$ (which are proper and smooth). Let $(NS^p)^{\bullet}$ and $(J^p_a)^{\bullet}$ denote the complexes induced by the simplicial structure and similarly to \eqref{NSext} we obtain an extension of $(NS^p)^{\bullet}$  by $(J^p_a)^{\bullet}$.  By taking homology groups we
then get boundary maps $$\lambda^i_a : H^i((NS^p)^{\bullet}) \to
H^{i+1}((J^p_a)^{\bullet}).$$ 
\subsubsection{\bf Hodge Conjecture (\protect{\cite[2.3.4]{BA}})}\label{HCconj} {\em The boundary map $\lambda^i_a$ behave well with respect to the extension class map $e^p$ yielding a {\em motivic} cycle class map, \ie the following diagram
$$\begin{array}{ccc}
H^{i}((NS^p)^{\bullet})&\by{\lambda^i_a}& H^{i+1}((J^p_a)^{\bullet})\\
\downarrow & &\downarrow\\
H^{2p+i}(X)^{p,p} & \by{e^p} & J^p(H^{2p+i}(X))
\end{array}$$
commutes.\footnote{Note that all maps in the square are canonically defined.} The image 1-motive (up to isogeny) is the Hodge 1-motive $\Xi^{i,p}(X)$ corresponding to $H^{2p+i}(X,\Z)^h$.\\ }

Moreover, one might then guess that the complex of 1-motivic sheaves $(CH_{\Xs}^p)^{(1)}$ would provide such Hodge 1-motive directly. 

\subsection{Non-homotopical invariant theories} A typical problem occurring with  homotopical invariant theories attached to singular varieties is that they do not catch some informations coming from the singularities. In general, the cohomological Picard 1-motive $\Pic^+ (X)$ of a proper scheme $X$ is given by the semi-abelian quotient of $\Pic^0 (X)$ (see {\it Scholium} \ref{surj}). Loosing its additive components we loose informations, \eg we don't see cusps. In order to reach the full picture here we have to enlarge our target to Laumon's 1-motives at least.  A natural guess is that
our 1-motives are only the \'etale part of Laumon's 1-motives, \ie there exists $\Pic^{+}_a(X,Y; i)\in \M^a$ such that 
$$\Pic^{+}_a(X,Y; i)_{\et} = \Pic ^{+}(X,Y; i)$$
and similarly $\RPic_a (X)\in D^b(\M^a)$ (\cf \ref{LRthm}), $\Xi^{i,p}_a(X)\in \M^a$ such that $\Xi^{i,p}_a(X)_{\et} = \Xi^{i,p}(X)$ (\cf \ref{AHCconj}), etc.\footnote{Note that such  Laumon 1-motives should rather be visible from a triangulated viewpoint! There should be a ``sharp" {\it cohomological}\, motive $M_{\sharp}(X)$ in a triangulated category $\DM_{\sharp}$, related to Voevodsky category of motivic complexes, with a realisation in $D^b({\rm FHS})$. The conjectural formalism for motivic complexes should be translated for $\sharp$-motivic complexes.} 

Their geometrical sources are {\it additive} Chow groups and their universal regular quotients, \cf \cite{BE} and, by the way, see \cite{ESV} for a construction of an additive version of the cohomological Albanese $\Alb^+_a (X)$ of a projective variety $X$, \ie here $\Alb^+_a (X)_{\et}= \Alb^+ (X)= \LA{1}^*(X)$, etc. as above. Similarly, for $X$ quasi-projective, we expect a  formal part defining $\Alb^+_a (X)$ as a Laumon 1-motive.

\subsubsection{ }  The forthcoming theories are {\it sharp}\,  cohomology theories, \eg $\sharp$-singular cohomology $X\leadsto H^*_{\sharp}(X)\in {\rm FHS}$ for $X$ over $\C$, $\sharp$-De Rham cohomology $H^*_{\sharp-DR}(X)$ of $X$ $k$-algebraic over a field $k$ of zero characteristic and $\sharp$-crystalline cohomology in positive characteristics, which are  non-homotopical invariant theories. Sample $$H^1_{\sharp}(X)\df T_{\oint}(\Pic^+_a (X))$$ Here $\Pic^+_a (X) = [0\to \Pic^0 (X)]$ if $X$ is proper: in this case, define the group scheme $\Pic^{\sharp}$ by the following pull-back square
 (\cf \eqref{sharp}, \ref{natural}, \ref{surj} and \cite[4.5]{BSAP}) 
\begin{equation}\label{sharpDR}
\begin{array}{ccccc}
 &\bPic^{\natural} (\Xs) & \to& \bPic (\Xs)&\\
&\uparrow &&\uparrow&\\
 &\Pic^{\sharp} (X) & \to& \Pic (X)&\\
\end{array}
 \end{equation}
 such that
 \begin{itemize}
 \item $\ker (\Pic^{\sharp, 0} (X) \onto \Pic^0 (X))=H^0(\Xs, \Omega^1_{\Xs})$ and
 \item $\ker (\Pic^{\sharp, 0} (X) \onto\bPic^{\natural, 0} (\Xs))$ is the additive subgroup  $\subseteq \Pic^0 (X)$; 
\end{itemize}
then 
 $$H^1_{\sharp-DR}(X)\df \Lie \Pic^{\sharp, 0} (X)$$
 so that $H^1_{\sharp-DR}(X)$ is an extension of $H^1_{DR}(X)$ by the additive part of $\Pic^0 (X)$.

\subsubsection{ } Similarly, remark that (see \cite{BSE}) for $E =  (H, V_{\d}) \in {\rm EHS}$ there is a surjection $$\Ext_{\rm EHS} (\Z (0), E) \onto \Ext_{\rm MHS}  (\Z(0), H)$$ and the kernel of this map is a vector space if $H = 
H^{2r-1}(X, \Z (r))$ and
$V_{i} = \HH^{2r-1}(X,\cO_X\to \cdots\to \Omega_X^{i -1}) (r)$, where $X$  is a proper $\C$-scheme; in particular, if $X$ is the cuspidal curve then $\Ext_{\rm EHS} (\Z (0), E)$ is the additive group $\G_a =\Pic^0 (X)$.

\subsection{Final remarks}
Hoping to have puzzled the reader enough to procede on these  matters I would finally remark that this exposition is far from being exhaustive. 

\subsubsection{ } For example,  for $S=\Spec(R)$ where $R$ is a complete discrete valuation ring and $K$ its function field,  a 1-motive over $K$ with good reduction (resp. potentially good reduction) is defined (in \cite{RAY}) by the property of yielding a 1-motive over $R$ (resp. after a finite extension of $K$). To any 1-motive $M=[L\to G]$ over $K$ is canonically associated (see \cite{RAY} for details) a {\em strict}\, 1-motive over $K$, \ie $M'=[L'\to G']$, such that $G'$ has potentially good reduction, a quasi-isomorphism $M'_{\rm rig}\to M_{\rm rig}$ in the derived category of bounded complexes of fppf-sheaves on the rigid site of  $\Spec (K)$, producing a canonical isomorphism  $T_{ \ell}(M')\cong T_{ \ell}(M)$ between the $\ell$-adic realizations, for any prime $\ell$. For a strict 1-motive $M=[L\to G]$ over $K$, the geometric monodromy  $\mu\colon  L\times T^{\vee}\to \Q$  (where $T^{\vee}$ is the character group of the torus $T\subseteq G$) is defined by valuating the trivialization of the Poincar\'e biextension. The geometric monodromy is zero if and only if $M$ has potentially good reduction. This theme is further investigated in \cite{BCV}. 

\subsubsection{ } The employ of 1-motives in arithmetical geometry is well testified, \eg see \cite{BR}, \cite{RI}, \cite{FC}, \cite{HS} and \cite{KT}. Note that in \cite{BMB} we also investigate $L$-functions with respect to Mordell-Weil and Tate-Shafarevich groups of 1-motives. Also the theme of 1-motivic Galois groups is afforded.  For $M$ a 1-motive over a field $k$ of zero characteristic let $M^{\otimes}$ be the Tannakian subcategory generated by $M$ in suitable mixed realisations (hopefully mixed motives).
The motivic Galois group of $M$, denoted $\Gal_{mot} (M)$, is the fundamental group of $M^{\otimes}$. The group $\Gal_{mot} (M)$ has an induced weight filtration $W_*$
and the unipotent radical $W_{-1} \Gal_{mot} (M)$ has a nice characterisation (see \cite{BG} for details). Furtehrmore, Fontaine's theory relating $p$-adic mixed Hodge structures over a finite extension $K$ of $\Q_p$ to mixed motives would provide categories of 1-motives over the $p$-adic field $K$ (see \cite{FJ}).  

\subsubsection{ } Passing from 1-motives to 2-motives is conceivable but (even conjecturally) harmless. A general guess is that there should be abelian categories  $$\mathcal{M}_0 \subseteq \M \subseteq \cdots \subseteq \mathcal{M} $$ where $\mathcal{M}_0=$ Artin motives, $\M$ = 1-motives and further on we have categories of $n$-motives $\mathcal{M}_n$  which can be realized as Serre subcategories of cohomological dimension $\leq n$ of the abelian category $\mathcal{M}$ of mixed motives. Assuming the existence of $\mathcal{M}$ a source of inspiration is \cite{GM}, \cite{V} and \cite{BEI}: such $\mathcal{M}_n$ would be somehow `generated' by motives of varieties of dimension $\leq n$ and $M (X)$, the  motive of $X$ smooth and projective, decomposes as $\oplus M^i(X)[-i]$ where $M^i(X)\in \mathcal{M}_i$ such that $M^i(X)= M^{2d-i}(X)$ for $d =\dim (X)$. 

\subsection*{Acknowledgement} I would like to thank all my co-authors, in particular: F. Andreatta, A. Bertapelle, M. Bertolini, J. Ayoub, B. Kahn and M. Saito for useful discussions and a keen interest in these matters.\footnote{Note that work in progress and preliminary versions of my papers are firstly published on the web and currently updated, \eg browsing from my home page.}

\begin{comment}
Here $R_{Hodge}(M^i (X)) = H^i(X)\in \mbox{MHS}$ compatibly as follows
$$\begin{array}{ccc}
\mathcal{M}_i(\C)&\longby{T}& \mbox{MHS}_i\\
\downarrow & &\downarrow\\
\mathcal{M}(\C)&\longby{R}& \mbox{MHS}
\end{array}$$
The wish is always that any realization of an algebraic variety which is of level $\leq i$ would be lifted back {\it via}\, $T$ and algebraically defined in $\mathcal{M}_i$. 
\end{comment}

\end{document}